\documentclass[12pt]{article}
\date{\empty}

\usepackage{graphicx}
\usepackage{multicol}
\usepackage{amsfonts}

\newcommand{\be}{\begin{equation}}
\newcommand{\ee}{\end{equation}}

\def\ds{\displaystyle}

\def \myfigures #1#2#3#4#5#6#7#8
{\begin{figure}[ht]
    \begin{center}
        \includegraphics[width=#2 \textwidth]{#1.eps}
        \includegraphics[width=#6 \textwidth]{#5.eps}
        \parbox[t]{#4\textwidth}{\caption {#3}\label{#1}}
        \parbox[t]{#8\textwidth}{\caption {#7}\label{#5}}
    \end{center}
\end{figure} }

\begin{document}
\begin{center}
\Large{\bf Structures and waves \\ in a  nonlinear heat-conducting medium}\\
%\maketitle
%%\author
~\\
\large{Stefka Dimova, \\
Faculty of Mathematics and Informatics, University of Sofia, Bulgaria\\
{\tt dimova@fmi.uni-sofia.bg}\\
%%\and
Milena Dimova, Daniela Vasileva, \\
Institute of Mathematics and Informatics, Bulgarian Acad. Sci.}\\
{\tt mkoleva,vasileva@math.bas.bg}
\end{center}
%%%%\institute{Stefka Dimova \at Faculty of Mathematics and Informatics, University of Sofia,
%%%%5 James Bourchier Blvd., 1164 Sofia, Bulgaria,
%%%%\email{dimova@fmi.uni-sofia.bg}
%%%%\and Milena Dimova \at Institute of Mathematics and Informatics, Bulgarian Acad. Sci.,
%%%%Acad. G.Bonchev Str., bl.8, 1113 Sofia, Bulgaria, \email{mkoleva@math.bas.bg}
%%%%\and Daniela Vasileva \at Institute of Mathematics and Informatics, Bulgarian Acad. Sci.,
%%%%Acad. G.Bonchev Str., bl.8, 1113 Sofia, Bulgaria, \email{vasileva@math.bas.bg}}

~\\~\\
\begin{tabular}{ll}
~~~~ & 
 \parbox[t]{0.86\textwidth}
{\large  The final publication will appear in Springer Proceedings in Mathematics and Statistics 
({\tt http://www.springer.com/series/10533}), Numerical Methods for PDEs:
Theory, Algorithms and their Applications.}
\end{tabular}
~\\~\\

\abstract{The paper is an overview of the main  contributions of a Bulgarian team of researchers to
the problem of finding the possible structures and waves in the open
nonlinear heat conducting medium, described by a reaction-diffusion
equation. Being posed and actively worked out by the Russian school of
A. A. Samarskii and S.P. Kurdyumov since the seventies of the last century, this problem
still contains open and challenging questions.}

\section{Introduction}
\label{sec:1}

A very general form of the model of heat structures reads as
follows:
\be
u_t = \sum_{i=1}^{N}(k_i(u)u_{x_i})_{x_i} + Q(u), \
t>0, \ x \in \mathbb{R}^N, \label{01}
\ee
where  the heat conductivity coefficients   $k_i(u) \ge 0$ and the heat source
$Q(u) \ge 0$ are nonlinear functions  of the temperature $u(t,x) \ge 0$.

Models such as  (\ref{01}) are  studied  by many researchers in various contexts.
A part of this research is devoted to semilinear equations:
$k_i(u)\equiv 1,$ $Q(u)=\lambda e^u$ (Frank-Kamenetskii equation) and
 $Q(u)=u^\beta,$ $\beta>1$. After the pioneer work of Fujita  \cite{79},
 these equations and some  generalizations of theirs are studied intensively
 by many authors, including  J. Bebernes,  A. Bressan, H. Brezis, D.
Eberly, A. Friedman, V.A. Galaktionov, I.M. Gelfand, M.A. Herrero,
R. Kohn, L.A. Lepin, S.A. Posashkov, A.A. Samarskii, J.L. V\'azquez, J.J.L. Vel\'azquez.
The book \cite{80} contains a part of these investigations and  a large bibliography.

The quasilinear equation is studied by D.G. Aronson,  A.
Friedman, H.A. Levine, S. Kaplan, L.A. Peletier, J.L. V\'azquez  and others. The contributions of the Russian school are
significant. The unusual localization effect of the blow-up boundary regimes is discovered
by numerical experiment in the work of A.A. Samarskii and M.I. Sobol in 1963
 \cite{22}. The problem of localization for quasilinear equations with a source is posed by S.P. Kurdyumov
\cite{13} in 1974. The works of I.M. Gelfand, A.S. Kalashnikov, the scientists of the school of A.A. Samarskii
and S.P. Kurdyumov are devoted to the challenging  physical and mathematical problems,
related with this model and its generalizations. Among them are:
localization  in space of the process of burning,
different types  of blow-up, arising of structures -
traveling and standing waves, complex structures  with varying degrees of symmetry.
The combination of the computational experiment with the progress
in the qualitative and analytical methods of the
theory of  ordinary and partial differential equations,  the Lie and the Lie-B\"{a}cklund
group theory, has been crucial for the success of these investigations.
The book \cite{32} contains many of these results, achieved to
1986, in the review  \cite{NEW} there are citations of later works.

A special part of these investigations is devoted to finding and studying
different kinds of self-similar and invariant
solutions of  equation (\ref{01}) with power nonlinearities :
\be k_i(u)=u^{\sigma_i}, \ \ Q(u)=u^ { \beta}. \label{02}
\ee

This choice is suggested by the following reasoning.

{\bf First}, such temperature dependencies are usual for many real processes  \cite{25}, \cite{86}, \cite{14}.
For example, when
$ { \sigma _i}= {\sigma}=2.5$, $ { \beta} \le 5.2$,
equation (\ref{01}) describes thermo-nuclear combustion in plasma in
the case of electron heat-conductivity; the parameters ${\sigma}=0$, $ \ 2 \le { \beta} \le 3 $
correspond to the models of autocatalytic processes with diffusion in the chemical reactors;
 $ { \sigma} { \approx} 6.5$ corresponds to the radiation heat-conductivity of the
 high-temperature plasma in the stars, and so on.

{\bf Second}, it is shown in \cite{27}, that in the class of power functions
the symmetry of  equation (\ref{01}) is maximal in some sense - the equation admits
a rich variety of invariant solutions. In general, almost all of the
dissipative structures known so far are invariant or partially invariant solutions of nonlinear equations.
The investigations of the dissipative structures
provide reasons to believe that the invariant solutions
describe the attractors of the dissipative structures' evolution and thus they characterize  important
internal properties of the nonlinear dissipative medium.

{\bf Third}, this rich set of invariant solutions of equation (\ref{01}) with power nonlinearities is
necessary for the successful application of the  methods for investigating the same equation
in the case of more general dependencies  $k_i(u)$, $Q(u)$. By using the methods of operator comparison
\cite{48} and stationary states  \cite{47} it is possible to analyze the properties of the solutions
(such as localization, blow-up, asymptotic behavior) of whole classes general nonlinear equations.
The method of approximate self-similar solutions
 \cite{31}, developed in the works of A.A. Samarskii and V.A. Galaktionov,
makes it possible to put in accordance with such general equations some other, basic equations.
The latter could have invariant solutions even if the original equations do not have such.
Moreover, the original equations may significantly differ from the basic equations,
and nevertheless their solutions
tend to the invariant solutions of the basic equations at the asymptotic stage.

{\bf Finally}, in the case of power coefficients the dissipation and the source are coordinated so that
complex structures arise, moreover, a spectrum of structures, burning consistently,
occurs.

Below we report about the  main contributions of the Bulgarian research team: S.N. Dimova,
\frame{M.S. Kaschiev}, M.G. Koleva, D.P. Vasileva, T.P. Chernogorova, to the
problem of finding the possible evolution patterns  in the
 heat-conducting medium, described by the reaction-diffusion
equation (\ref{01}), (\ref{02}). The outline of the paper is as follows.
The main notions, needed further, are introduced in Section 2 on the simplest
and the most studied radially symmetric case. The specific peculiarities
of the numerical methods, developed and applied
 to solve the described problems, are systematized in Section 3. A brief report on the main
achievements of the team is made in Section 4. Section 5 contains some  open problems.

\section{The radially symmetric case, the main notions}
\label{sec:2}

Let us introduce the main notions to be used further on the Cauchy problem for equation
(\ref{01}) with initial data:
$$u(0,x)=u_0(x)\ge 0, \ x \in \mathbb{R}^N, \ \sup u_0(x) < \infty. $$
This problem could have global or blow-up solutions. The global in time solution  is
defined and bounded in $\mathbb{R}^N$ for every  t. {\it The unbounded
(blow-up) solution} is defined in $\mathbb{R}^N$ on a finite interval
$[0,T_0)$, moreover
$$\overline{\lim}_{t \to T_0^-} \sup_{x \in \mathbb{R}^N} u(t,x) = +\infty.$$
The time $T_0$ is called {\it blow-up time}.

The unbounded solution of the Cauchy problem with finite support initial data
$u_0(x)$ is called  {\it localized  (in a strong sense)}, if the set
$$
\Omega_L=\{x \in \mathbb{R}^N: \  u(T_0^-,x):=
\ds\overline{\lim}_{t\to T_0^-}  \; u(t,x)>0 \}
$$
is bounded in $\mathbb{R}^N$. The set $ \Omega_L$  is called  {\it localization region}.
The  solution localized in a strong sense grows infinitely for $t \to
T_0^-$ in a finite  region
$$\omega_L=\{x \in \mathbb{R}^N: \ u(T_0^-,x)= \infty \} $$
in general different from $\Omega_L$.

If for  $k_i(u) \equiv 0$ the condition
\be
\int_{1}^{\infty}{du \over Q(u)} < +\infty
\label{05}
\ee
holds, then the solution of the Cauchy problem is unbounded \cite{32}.
The heating of the medium happens in a blow-up regime, moreover the blow-up time
of every  point of the medium is different, depending on its initial temperature.

If for  $Q(u) \equiv 0$ the condition
\be
\int_{0}^{1} {k_i(u) \over u}\,du < +\infty, \ \
i=1,2,...,N, \label{06}
\ee
holds, then a {\bf finite speed} of heat propagation takes place for a finite support
initial perturbation in an absolutely cold medium \cite{32}.

In the case (\ref{02}) of power nonlinearities it is sufficient to have
 $ \sigma_i >0, \beta>1$ for  the conditions (\ref{05}) and (\ref{06}) to be satisfied. Then $k_i(0)=0$ and equation
(\ref{01}) degenerates. In general it has a generalized solution, which could have
discontinuous derivatives on the surface of degeneration $\{u=0\}.$

\bigskip
\noindent
{\bf 2.1 The basic blow-up regimes} will be explained on the radially
symmetric version of the Cauchy problem for equation (\ref{01}):
\be
u_t={1\over x^{N-1}} (x^{N-1}
u^{\sigma}u_x)_x + u^{\beta}, \quad  x\in \mathbb{R}^1_+, \quad t>0, \quad
\sigma>0, \ \beta >1, \label{07}
\ee
\be
u_x(t,0)=0, \quad u(0,x)=u_0(x)\ge 0,\ 0\le x<l, \quad u_0(x)\equiv 0, \quad
x\ge l.  \label{08}
\ee
If  $u_0(x)$ satisfies the additional conditions
$  u_0(x) \in C(\mathbb{R}_+^1),\ (u_0^\sigma u_0')(0)=0,$
there exists unique local (in time) generalized solution  $u=u(t,x)$ of problem
(\ref{07})-(\ref{08}), which is a nonnegative continuous function in $\mathbb{R}_+^1 \times (0,T),$
where $T \in (0,\infty]$  is the finite or infinite time of existence of the solution
(see the bibliography in the review \cite{130}). Moreover $u(t,x)$ is a classical solution
in a vicinity of every point $(t,x)$, where $u(t,x)$ is strictly positive. It could not have
the necessary smoothness  at the  points of degeneracy, but the heat flux  $-x^{N-1} u^{\sigma}u_x$
must be continuous. It means  that $u^{\sigma}u_x=0$ everywhere $u=0$.

Equation  (\ref{07}) admits a {\it self-similar solution (s.-s.s.)}  \cite{32}:
\be
u_s(t,x)=\varphi(t)\theta_s(\xi)=\left(1-{t\over
T_0}\right)^{-1\over\beta-1} \theta_s(\xi),\label{09}
\ee
\be
\xi=x/\psi(t)=x/\left(1-{t\over T_0}\right)^{m\over{\beta-1}}, \ \
m={\beta-\sigma-1\over 2}. \label{010}
\ee
The s.-s.s. corresponds to initial data $u_s(0,x)= \theta_s(x)$.
The function $\varphi(t)$ determines the amplitude of the solution.
The {\it self-similar function (s.-s.f.)} $\theta_s(\xi)\ge 0$ determines
the space-time structure of the s.-s.s. (\ref{09}).
This function satisfies the degenerate ordinary differential equation in $\mathbb{R}_+^1$:
\be \label{011}
L(\theta_s) \equiv -{1\over \xi^{N-1}} (\xi^{N-1} \theta^{\sigma}_s
\theta'_s)' +{\beta-\sigma-1\over 2(\beta-1)T_0} \xi\theta'_s
+{1\over (\beta-1)T_0}\theta_s-\theta^{\beta}_s=0  \label{13}
\ee
and the boundary conditions:
\be\label{012}
\theta'_s(0)=0,\ \  \theta_s(\infty)=0,\ \ \theta^{\sigma}_s\theta'_s(\xi_0)=0,\
\mbox{if}\ \ \theta_s(\xi_0)=0.
\ee

Equation (\ref{011}) has two constant solutions:
$\theta_s(\xi)\equiv \theta_H =$  $(T_0(\beta-1))^{-1\over \beta-1}$ and $\theta_s(\xi)\equiv 0.$
These two solutions play an important role in the analysis of the different solutions of equation (\ref{011}).
For blow-up regimes we assume $T_0>0.$ Without loss of generality we set
\be
T_0= 1/ (\beta-1), \ \ \mbox{then} \ \ \theta_H\equiv 1. \label{014}
\ee

The analysis of the solutions of problem (\ref{011}), (\ref{012}), carried out in the works
\cite{36}, \cite{10},  \cite{39}, \cite{41}, (see also \cite{32}, Chapter IV),
gives the following results:
\begin{itemize}

\item
For arbitrary $1<\beta\le \sigma+1$ there exist a finite support solution
$\theta_s(\xi)\ge 0.$

\item For $\beta<\sigma +1, \ N\ge1$ and $\beta =
\sigma +1, \ N>1$ the problem has no nonmonotone solutions.
The uniqueness is proved only for $\beta<\sigma+1,\ N=1.$

The graphs of the s.-s.f. $\theta_s(\xi)$ for $\beta=\sigma+1=3, \ N=1,2,3$ are shown in Fig. 1,
the graphs of the s.-s.f. $\theta_s(\xi)$ for
$\beta=2.4<\sigma+1=3, \ N=1,2,3$ -- in Fig. 2.

%%\myfigures{s2b3n}{0.45}{\footnotesize{$\beta = \sigma +1$,
%%$S$-regime}}{0.45} {s2b4n}{0.45}{\footnotesize{$\beta < \sigma +1$,
%%$HS$-regime }}{0.45}

\myfigures{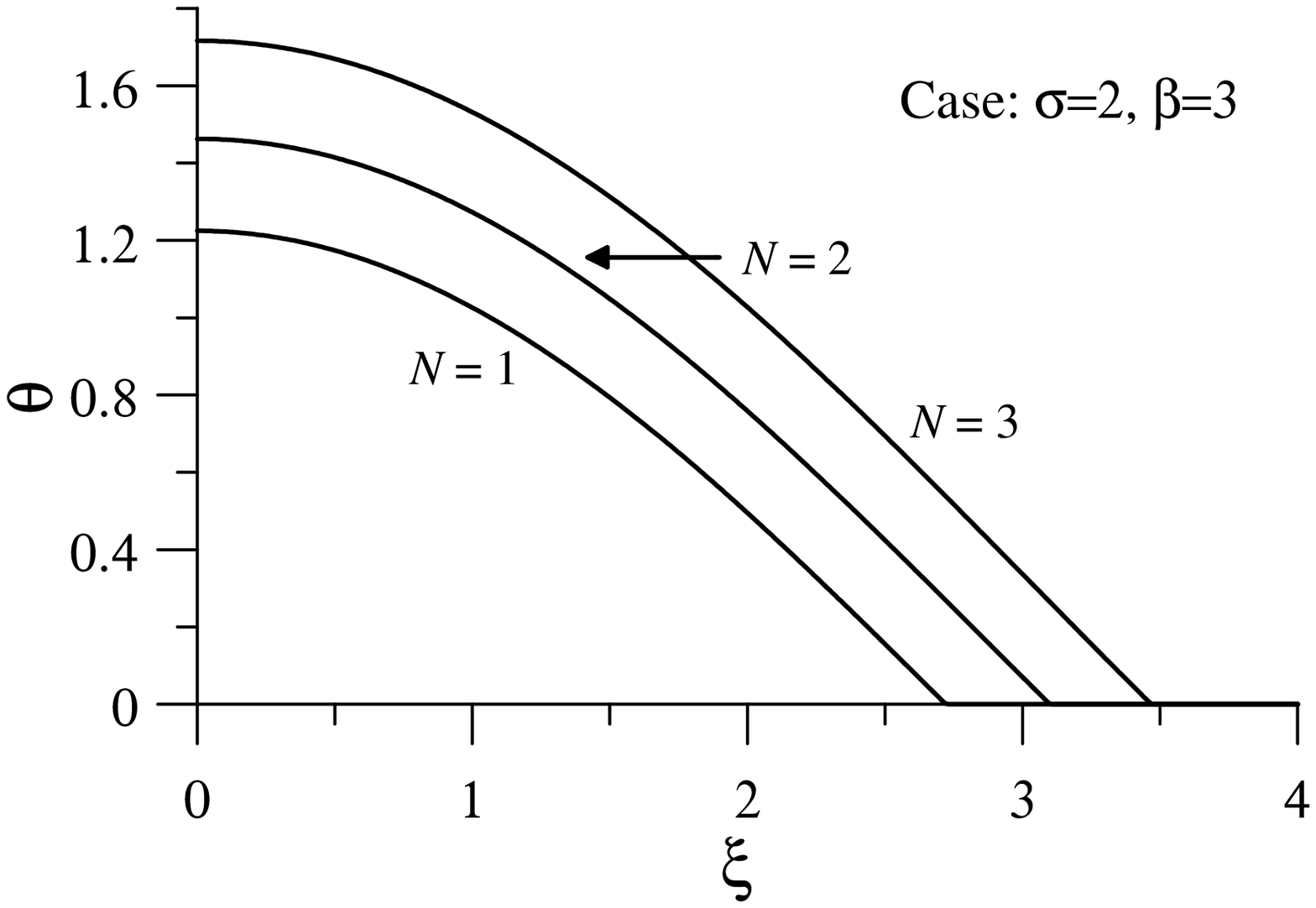}{0.4}{$\beta = \sigma +1$,
$S$-regime}{0.4} {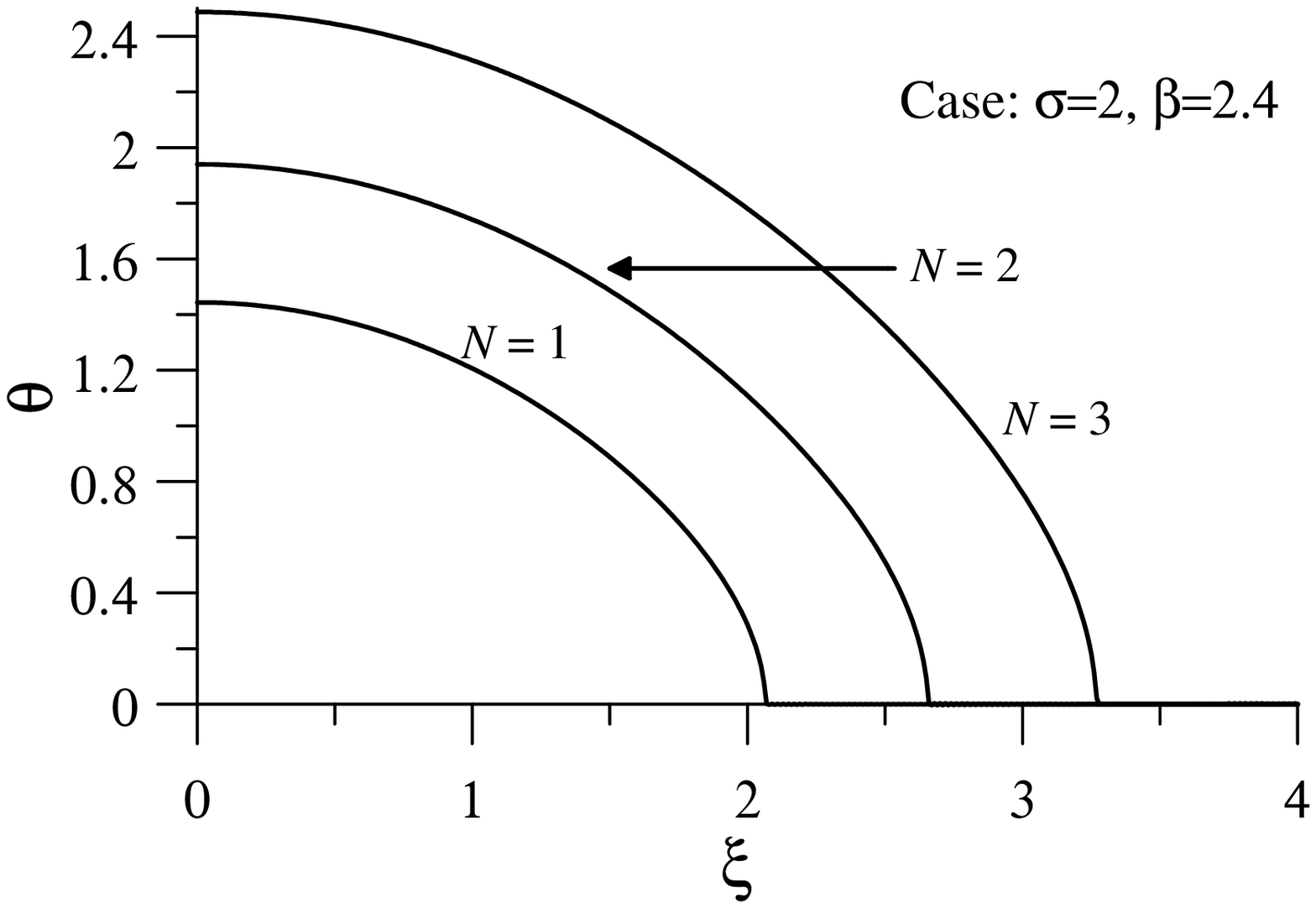}{0.4}{$\beta < \sigma +1$,
$HS$-regime }{0.4}

\item For $\beta>\sigma+1, \ N\ge 1$ the problem has no finite support solutions.

\item If $\sigma+1 < \beta < \beta_s=(\sigma+1)(N+2)/{(N-2)_+}$,
{\it ($\beta_s$ -- the critical Sobolev exponent)}, the problem has at least one solution
$\theta_s(\xi)>0$ in $\mathbb{R}^1_+$, strictly monotone decreasing in
 $\xi$ and having the asymptotics
\be
\theta_s(\xi)=C_s \xi^{-2/(\beta-\sigma-1)}[1+\omega(\xi)],
\qquad \omega(\xi)\to 0, \qquad \xi\to\infty,\label{015}
\ee
$C_s=C_s(\sigma,\beta,N)$ is a constant. Later on in \cite{50} the interval in $\beta$ has been extended.

\item For $N=1,\ \ \beta>\sigma+1$ the problem has at least
\be
K=-[-a]-1,\ \ a={\beta-1\over \beta-\sigma-1}>1\label{016}
\ee
different solutions (\cite{41}, \cite{39}, \cite{32}). Let us introduce the notations $\theta_{s,i}
(\xi),\ i=1,2, \ldots , K$ for them. On the basis of linear analysis   and some numerical results
in the works \cite{43}, \cite{57} it has been supposed that the number of different solutions  $ \theta_{s,i} (\xi)$
for  $\beta>\sigma +1$ and  $N\ge 1$ is $K+1.$  For $N=1$ this result was refined \cite{132}
by using bifurcation analysis: the number of solutions is
$K=[a]$, if $a$ is not an integer, and $K=a-1$, if $a$ is an integer.
For $N=2,3$ the bifurcation analysis gives the same estimate for the number of different solutions,
but for $\beta \approx \sigma+1$, $\beta> \sigma+1$ it is violated (see 4.1).
\end{itemize}
\myfigures{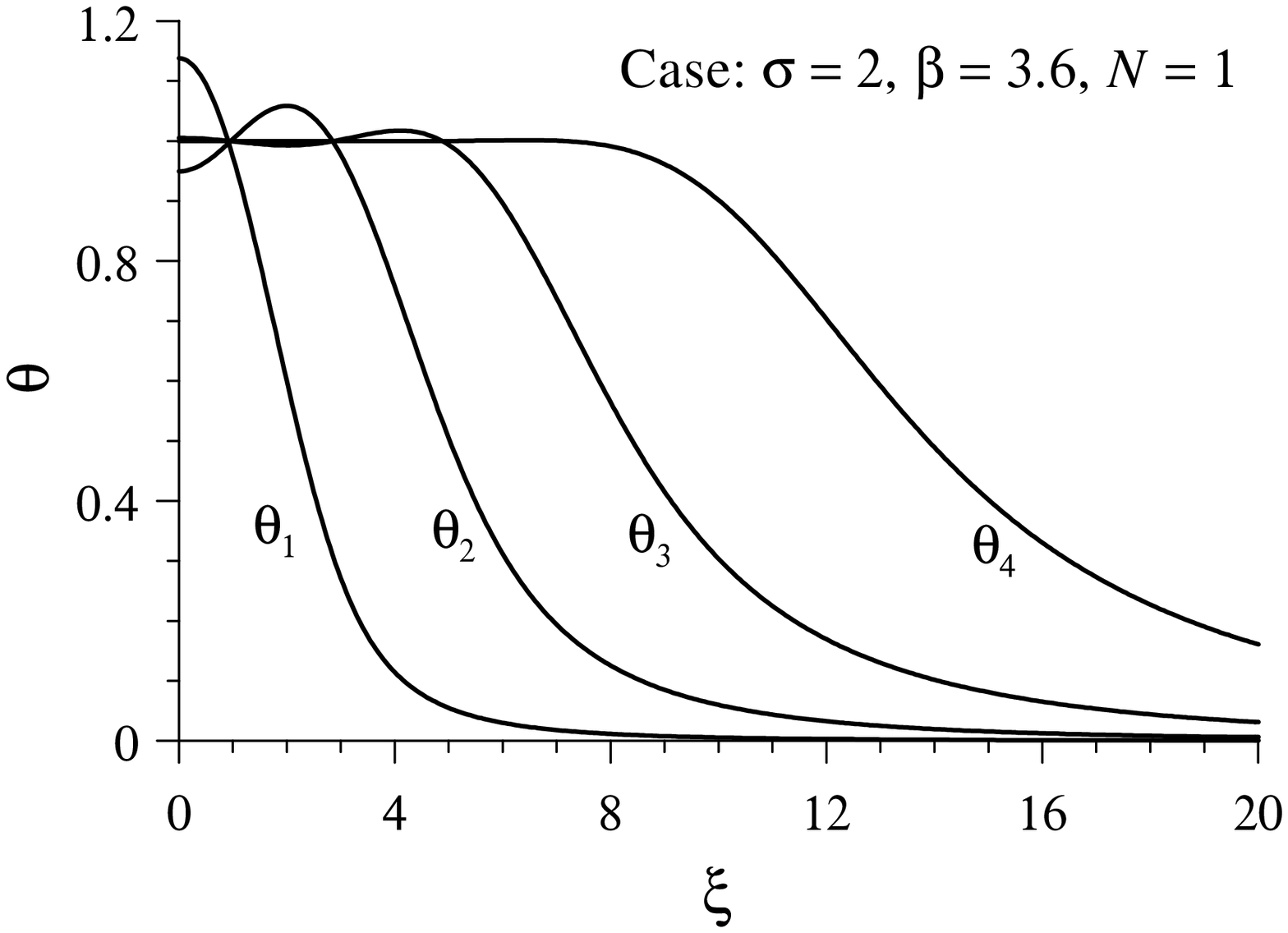}{0.4}{$\beta > \sigma +1$, $LS$-regime}{0.4}
{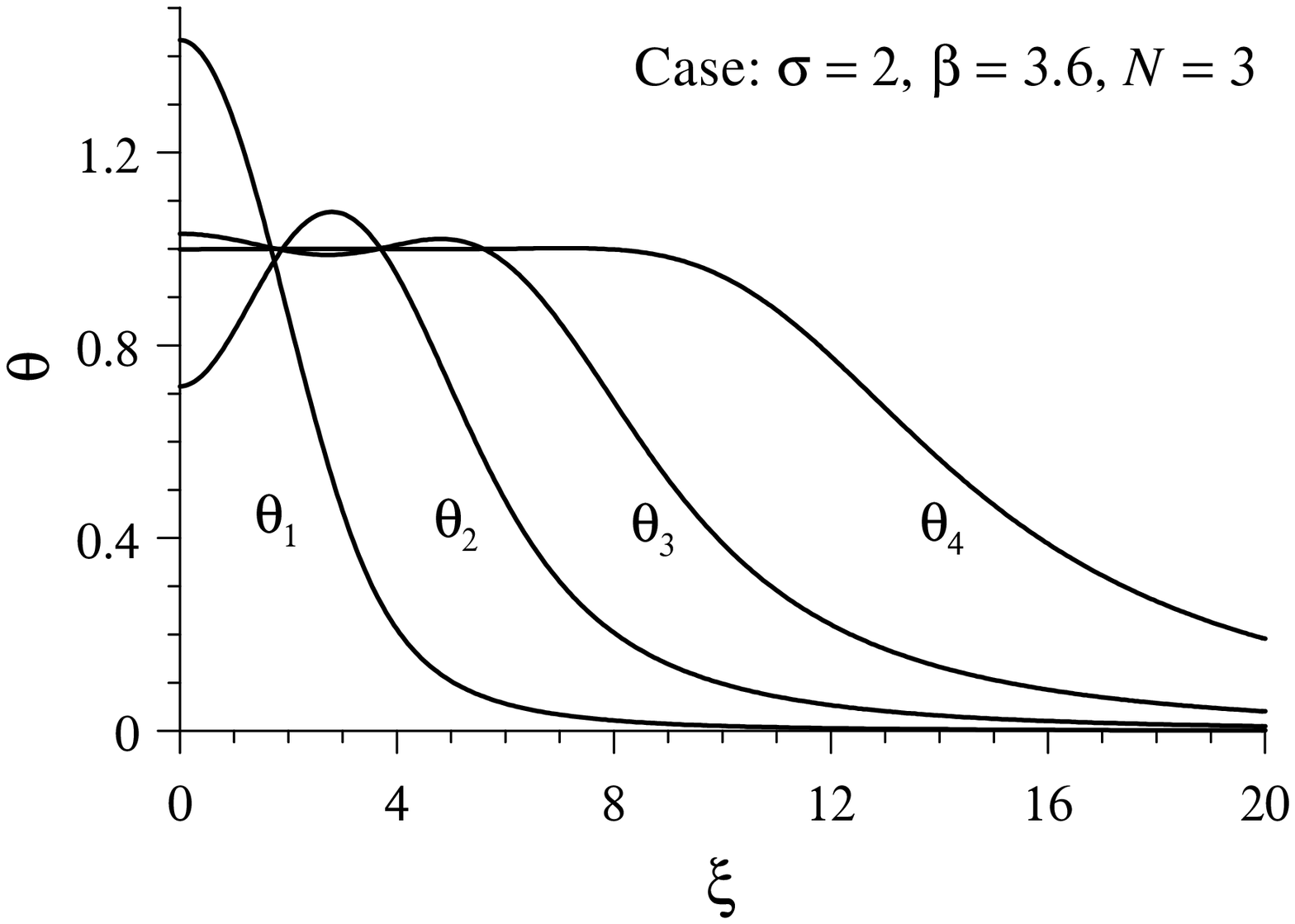}{0.4}{$\beta > \sigma +1$, $LS$-regime}{0.4}

The graphs of the four self-similar functions, existing for
  $\sigma=2, \beta=3.6>\sigma+1$  ($K=4$) are shown in  Fig. 3 ($N=1$) and Fig. 4 ($N=3$).

These results determine {\bf the basic regimes of burning of the medium},
described by the s.-s.s.  (\ref{09}), (\ref{010}). The following notions
are useful for their characterization:

-- {\it semi-width} $x_s=x_s(t)$, determined by the equation
$u(t,x_s)=u(t,0)/2$ for  solutions, monotone in $x$ and having
a single maximum  at the point $x=0$;

-- {\it front-point} $x_f$: $u(t,x_f)=0,\ \ u^{\sigma} u_x(t,x_f)=0.$

\medskip
\noindent
{\bf 2.1.1 $HS$-evolution, total blow-up, $1<\beta<\sigma+1$}

\noindent
The heat diffusion is more intensive than
the heat source. The semi-width  and the front tend to infinity; a {\it heat wave},
which covers the whole space for time
 $T_0,$  is formed. The process is not localized:
$ \displaystyle {\rm mes}\ \Omega_L = {\rm mes}\ \omega_L=\infty, \ \ x_s\to \infty, \ \
x_f\to \infty, \ \ t\to T^-_0. $

\medskip

\noindent
{\bf 2.1.2 $S$-evolution, regional blow-up, $\beta=\sigma+1$}

\noindent
The heat diffusion and the  source are correlated in such a way,
that leads to localization of the process in a region
$\Omega_L = \omega_L =\left\{ |x|<{L_s / 2}\right\}$
of diameter  $L_s,$ called a {\it fundamental length} of the $S$-regime.
The semi-width is constant, inside $\Omega_L$ the medium is heated to infinite temperature for time
 $T_0:\ x_s={\rm const} , \ x_f=L_s/2.$
In the case  $N=1$ the solution $\theta_s(\xi)$ ({\it Zmitrenko-Kurdyumov solution}) is found  \cite{10} explicitly:
\be
\theta_s(\xi)=\left\{
\begin{array}{ll}
\displaystyle\left( {2(\sigma+1)\over\sigma+2} \cos^2 {\pi\xi\over
L_s}\right)^{1/\sigma}, &\displaystyle |\xi|\le {L_s\over 2}\\
0 &\displaystyle |\xi|>{L_s\over 2},\end{array}\right.
\label{019}
\ee
$L_s={\rm diam} \ \Omega_L=(2\pi \sqrt{\sigma+1}) / \sigma,\ \ x_s=L_s
\arccos ((2^{-{\sigma\over 2}})/\pi).$

The solution (\ref{019}) is called
{\it elementary solution} of the $S$-regime for $N=1$. In this case equation (\ref{011}) is autonomous and every function,
consisting of  $k$ elementary solutions, $k=1,2,\ldots $, is a solution as well, i.e.,
equation (\ref{011}) has a countable set of solutions.

\medskip

\noindent
{\bf 2.1.3 $LS$-evolution, single point blow-up, $\sigma+1<\beta<\beta_f=\sigma+1+{2\over N}$}

\noindent
Here {\it $\beta_f$ is the critical Fujita exponent} \cite{81}.
The intensity of the source is bigger, than the diffusion.
The front of the s.-s.s. is at infinity
 (\ref{015}), the semi-width decreases and the medium is heated to infinite temperature in a single point:
$${\rm mes}\  \omega_L=0,\ \ x_s \to 0,\ \ t\to T^-_0.$$
According to the different s.-s.f.  $\theta_{s,i}(\xi),\
i=1,2,\ldots,$ the medium burns as a {\it simple} structure $(i=1)$ and as {\it complex} structures
$(i>1)$ with the same blow-up time.

\medskip

\noindent
{\bf 2.2 Stability of the self-similar  solutions}

\noindent
To show the important property
of the s.-s.s. as attractors of wide classes of other solutions of the same
equation,  we will need of  additional notions.

In the case of arbitrary  finite support initial data
$u_0(x)$ (\ref{08}) the so called  {\it self-similar representation} \cite{39} of the solution
 $u(t,x)$ of problem (\ref{07}), (\ref{08}) is defined.
 It is determined at every time  $t$ according to the structure of the s.-s.s.
(\ref{09}), (\ref{010}):
\be
\Theta(t,\xi)  = \left(1-{t /
T_0}\right)^{1\over \beta-1} u\left(t,\xi\left(1-{t / T_0}\right)^{m \over \beta-1}\right)
 =\varphi^{-1}(t)u(t,\xi \psi(t)). \label{020}
\ee
The s.-s.s. $u_s(t,x)$ is called  {\it asymptotically stable} \cite{32},
if there exists a sufficiently large class of solutions  $u(t,x)$ of problem (\ref{07}),
(\ref{08}) for initial data  $u_0(x) \not \equiv \theta_s (x)$, whose self-similar representations
  $\Theta(t,\xi)$ tend in some norm to $\theta_s(\xi)$ when $t\to T_0^-$:
\be
\|\Theta(t,\xi)-\theta_s(\xi)\|\to 0, \ t\to T^-_0. \label{021}
\ee

 The definition of the self-similar representation (\ref{020}) contains the blow-up time $T_0.$
For theoretical investigations this is natural, but for numerical investigations definition (\ref{020}) is unusable
since for arbitrary initial data  $u_0(x) \neq \theta_s(x)$ $T_0$ is not known.
Therefor another approach has been proposed and numerically implemented  (for $N=1$) in \cite{36}, \cite{39}.
This approach  gives a possibility to investigate {\bf the structural stability} of the unbounded solutions
in a special {\bf``self-similar'' norm},  consistent for every $t$ with
 the geometric form of the solution and not using explicitly the blow-up time $T_0$.
A new self-similar representation, consistent with the structure of the s.-s.s. (\ref{09}), (\ref{010}) is introduced:
\be\label{024}
\Theta(t,\xi)=u(t,\xi(\gamma(t))^{-m})/\gamma(t),
\qquad \gamma (t)={{\max_x} u(t,x)\over {\max_ \xi} \theta_s(\xi)}.
\ee
If the limit (\ref{021})
takes place for  $\Theta(t,\xi)$, given in (\ref{024}), then the self-similar solution
$u_s(t,x)$ is called  {\it structurally stable}.

The notion of structural stability, i.e., the preservation in time
of some characteristics of the structures, such as geometric form, rate of growth, localization in space,
is tightly connected with the notion invariance of the solutions with respect to the transformations,
involving the time \cite{33}. This determines its advisability for investigating
the asymptotic behavior of the blow-up solutions.

In the case of complex structures another notion of stability is needed, namely {\bf metastability}.
 The self-similar solution  $u_s(t,x)$ is called  {\it metastable},
if for every  $\varepsilon>0$ there exists a class of initial data
$u(0,x)\approx \theta_s(x)$ and a time $T$, $T_0-T\ll T_0$ such that
$$\parallel\Theta(t,\xi)-\theta_s(\xi)\parallel
\le\varepsilon , \ \ \mbox{for}\ \ 0\le t\le T$$
holds for the self-similar representations  (\ref{024}) of the corresponding solutions.
This means, that the metastable s.-s.s. preserves its complex space-time structure during the evolution up to time
 $T$, very close to the blow-up time
 $T_0$. After that time the complex structure could degenerate into one or several simple structures.

\section{Numerical methods}
\label{sec:3}
To solve  the reaction-diffusion problem (\ref{07}), (\ref{08}) and
 the corresponding self-similar
problem  (\ref{011}), (\ref{012}), as well as their generalizations
both for systems of such equations and for the 2D case, appropriate
numerical methods and algorithms were developed.

The difficulties, common for the nonstationary and for the self-similar
problems were: the nonlinearity, the dependence on a number of
parameters (not less than 3); insufficient smoothness of the
solutions on the degeneration surface, where the solutions vanish. In
the case of radial symmetry for $N>1$ and polar coordinates in the
2D case additional singularity at $x=0 \ (\xi=0)$ occurs.

The main challenge  in solving the self-similar problems is the
non-uniqueness of  their solutions for some ranges of the parameters. The following problems arise: to
find a ``good'' approximation to each of the solutions; to
construct an iteration process, converging fast to the desired  solution
(corresponding to the initial approximation) and
ensuring sufficient accuracy; to construct a computational process,
which enables  finding all different solutions
for given parameters $(\sigma, \beta, N)$ in one and the same way; to
determine in advance where to translate the boundary conditions from
infinity, for example (\ref{012}), for the asymptotics  (\ref{015})
to be fulfilled.

The main difficulty in solving the nonstationary problems is the
blow-up  of their solutions:  blow-up in a single
point, in a finite region and in the whole space. And two others,
related with it: the moving front of the solution, where it is often
not sufficiently  smooth; the instability  of the blow-up solutions.

\medskip

\noindent
{\bf 3.1 Initial approximations to the different s.-s.f. for a given set of parameters}

\noindent
To overcome the difficulty with the initial approximations to
the different s.-s.f., we have used the approach, proposed and used in the works
 \cite{43},  \cite{57}. Based on the hypothesis that in the region of their
nonmonotonicity the s.-s.f. have small oscillations around the homogenious solution
  $\theta_H,$  this approach consists of  ``linearization''
of the self-similar equation around
  $\theta_H$ and followed by ``sewing'' the solutions of the resulting linear equation with
  the known asymptotics at infinity, e.g. (\ref{015}).

Our experiments showed \cite{D09}, that when $\beta\to \sigma+1+0$
the hypothesis about small oscillations of the s.-s.f. around $\theta_H$,
 is not fulfilled.
The detailed analytical and numerical investigations
\cite{D09}, \cite{D11}, \cite{D25} of the ``linear approximations''
in the radially symmetric case for
$N=1,2,3,$ showed that even in the case these approximations take
negative values in vicinity of the origin,
they still give the true number of  crossings with
 $\theta_H$ and thus, the character of nonmonotonicity of the s.-s. functions.
Recommendations of how to use the linear approximations in these cases are made in \cite{D09}.

Let us note, the ``linear approximations'' are expressed by different special functions:  the confluent
hypergeometric function $_1F_1(a,b;z)$ and the Bessel function $J_k(z)$  for different parameter
ranges    within the complex plane for $a$ and
$b$ and different ranges of the variable $z$.
To compute these special functions, various methods were used: Taylor series expansions,
expansions in ascending series of Chebyshev polynomials, rational approximations,
asymptotic series \cite{D16}, \cite{D24}.

\medskip

\noindent
{\bf 3.2 Numerical method for the self-similar problems}

\noindent
To solve the self-similar problem (\ref{011}), (\ref{012}) and its
generalization for systems of ODE and for the 2D-case,
the Continuous Analog of the Newton's Method (CANM) was used
(\cite{D25}-\cite{D09}, \cite{D11}, \cite{136}, \cite{D14}).
 Proposed by Gavurin  \cite{59}, this method was further developed  in \cite{60},
\cite{62} and used for solving many nonlinear problems.
The idea behind it is to reduce the stationary  problem
$L(\theta)=0$  to the evolution one:
\be\label{138}
L'(\theta)\frac {\partial \theta}{\partial t} = -L(\theta), \quad \theta(\xi,0)=\theta_0(\xi),
\ee
by introducing a continuous parameter $t$,  $0<t<\infty$,
on which the unknown solution depends: $\theta=\theta(\xi,t)$.
By setting $v=\partial\theta/\partial t$ and applying the Euler's method to the Cauchy problem
(\ref{138}),  one comes to the iteration scheme:
\be
L'(\theta_n)v_n=-L(\theta_n), \label{140}
\ee
\be
\begin{array}{l}
\theta_{n+1}=\theta_n+\tau_nv_n, \quad 0<\tau_n\le 1, \quad n=0,1,
\dots, \\[3pt]
\theta_n=\theta_n(\xi)=\theta(\xi,t_n), \quad v_n=v_n(\xi)=v(\xi,t_n), \\[3pt]
\qquad \theta_0(\xi) \ \mbox{being the initial approximation}.
\end{array}
\label{141}
\ee
The linear equations (\ref{140}) (or the system of such equations in the case of
a two-component medium)  are solved by the Galerkin Finite Element Method (GFEM) at every iteration step.
The combination of the CANM and the GFEM turned out to be very successful.
The linear system of the FEM with nonsymmetric matrix is solved by using
the  $LU$ decomposition. The iteration process (\ref{141}) converges very fast -- usually less than 15-16 iterations
 are sufficient for stop-criterium $ \|L(\theta_n)\| < 10^{-7}$.
The numerical investigation of the accuracy of the  method being implemented
shows errors (i) of order $O(h^4)$ when using quadratic elements in the radially symmetric case,
 and (ii) of optimal order $O(h^2)$ when using linear elements in the same case or bilinear ones in the 2D case.
 To achieve the same accuracy in vicinity of the origin
in the radially symmetric case for $N\ge3$ the nonsymmetric Galerkin method
 \cite{84} was developed  \cite{D11}, \cite{D23}.

The computing of the solutions of the linearized self-similar equation and their sewing
with the known asymptotics is implemented  in a software, so the process is fully automatized.
The software enables the computing of the self-similar functions for all of the blow-up regimes,
moreover in the case of $LS$-regime only the  number $k$ of the self-similar function $\theta_{s,k}(\xi)$ must be given.

\medskip

\noindent
{\bf 3.3 Numerical method for the reaction-diffusion problems}

\noindent
The Galerkin Finite Element Method (GFEM), based on
 the Kirchhoff transformation of the nonlinear heat-conductivity coefficient:
\be \label{g1}
G(u) = \int_{0}^{u}{s^\sigma}\,ds = u^{\sigma+1}/(\sigma+1),
\ee
was used (\cite{D01}, \cite{D03}, \cite{D16}, \cite{D25}-\cite{D26}, \cite{D18}-\cite{D21}, \cite{D19}, \cite{78})
for solving the reaction-diffusion problems.
This transformation is crucial for the further interpolation of the nonlinear coefficients
on the basis of the finite element space and for  optimizing
of the computational process.

Here below we point out the main steps of the method on the problem (\ref{07}), (\ref{08})
in a finite interval [0, $X(t)$] under the boundary condition $u(X(t))=0$.
Because of the finite speed of heat propagation we choose $X=X(t)$ so as to avoid the
 influence of the boundary condition on the solution.
The discretization is made on the Galerkin form of the problem:

Find a function $u(t,x)\in D,$
$$ D=\{u: x^{(N-1)/2}u, \ x^{(N-1)/2} {\partial
u^{(\sigma+1)/2}/ \partial x }\in L_2, \ u(X(t))=0 \},$$
which for every fixed $t$ satisfies the integral identity
\be\label{223}
(u_t,v)=A(t;u,v), \quad \forall v\in   H^1 (0,X(t)),\ \ \ 0<t<T_0,
\ee
and the initial condition (\ref{08}).

Here
$$
(u,v)=\int^{X(t)}_0 x^{N-1} u(x) v(x)dx, \ \ \ A(t;u,v)=\int^{X(t)}_0
\left[x^{N-1} \frac{\partial G(u)}{\partial x} \frac{\partial v} {\partial x}
+ x u^\beta v\right]dx,
$$
$$
 H^1 (0,X(t))=\{v:x^{(N-1)/2} v, \ x^{(N-1)/2} v'\in L_2(0,X(t)), \ v(X(t))=0. \}
$$
The lumped mass finite element method \cite{97} with interpolation
of the nonlinear coefficients $G(u)$ (\ref{g1}) and $q(u)=u^\beta$:
$$G(u) \sim G_I=\sum^n_{i=1} G(u_i)\varphi_i(x), \qquad  q(u) \sim q_I=\sum^n_{i=1} q(u_i)\varphi_i(x)$$
on the basis $\{\varphi_i\}, i=1,\dots,n,$ of the finite element space is used for discretization of (\ref{223}).
The resulting system of ordinary differential equations
with respect to the vector $U(t)=(u_1(t), u_2(t),\ldots,u_n(t))^T$ of the nodal values
of the solution  $u(t,x)$ at time $t$ is:
\be\label{2301}
\dot U= \tilde M^{-1}(-K G(U)) + q(U), \quad  U(0)=U_0.
\ee
Here the following denotations are used:
$ G(U)=(G(u_1),\ldots,G(u_n))^T, \quad  q(U)=(q(u_1),\ldots,q(u_n))^T,$
$\tilde M$ is the lumped mass matrix, $K$ is the stiffness matrix.
Let us mention, thanks to the Kirchhoff transformation and the interpolation
of the nonlinear coefficients only the two vectors $G(U)$ and $q(U)$ contain the
nonlinearity of the problem, while the matrix $K$  does not depend on the unknown solution.

To solve the system (\ref{2301}) an explicit
Runge-Kutta method \cite{66} of second order of accuracy and an extended region of
stability was used. A special algorithm for choosing the time-step $\tau$ ensures
the validity of the weak maximum principle  and,
 in the case of smooth solutions, the achievement of a given accuracy
 $\epsilon$ up to  the end of the time interval. The stop criterion is  $\tau<10^{-16}$ and then
$\tilde{T}_0$ is the approximate blow-up time, found in the computations.

It is worth mentioning, that the nonlinearity has changed the prevailing opinion about the explicit methods.
Indeed, there are at least two reasons an explicit method to be preferred over the implicit one
for solving the system (\ref{2301}):

-- the condition for solvability of the nonlinear discrete system on the upper time level imposes the same
restriction on the relation "time step -- step in space", as does the condition for validity of the weak maximum principle
for the explicit scheme (see \cite{32}, chapter VII, \S5);

-- the explicit method for solving large discrete systems has a significant
advantage over  the implicit one with respect to the computational complexity.

Let us also mention, that in the case of blow-up solutions the discrete system on the upper time level
would connect solution values differing by 6-12 orders of magnitude, which causes additional difficulties to overcome.
Finally the explicit methods allow easy parallelization.

 The special achievement of
 the proposed methods are the {\bf adaptive meshes} in the  $LS$-regime (refinement of the mesh)
 and in the  $HS$-regime (stretching meshes with constant number of mesh-points), consistent with the
self-similar low. Let us briefly describe this adaptation idea on the differential problem
\be
u_t= Lu, \ \ u=u(t,x), \ \ x \in \mathbb{R^N}, \ \ t>0,  \label{111}
\ee
which admits a self-similar solution of the kind
\be
u_s(t,x) = \varphi(t) \theta_s (\xi), \qquad \xi =  {x} / {\psi (t)}, \ \ u_s(0,x)=\theta_s(x).  \label{112}
\ee
Since the invariant solution $u_s(x,t) $
is an attractor of the solutions of equation (\ref{111}) for large
classes of initial data different from $\theta_s(x)$, it is
important to incorporate the structure (\ref{112}) in the
numerical method for solving equation (\ref{111}). The relation (\ref{112}) between $\xi$ and $x$
 gives the idea how to adapt the mesh in space. Let $\Delta x^{(0)}$ be the initial step in space,
$\Delta x^{(k)}$ -- the step in space at $t=t^k$. Then
 $\Delta x^{(k)}$ must be chosen so that $ \Delta \xi^{(k)}$
is bounded from below and from above
$$ { \Delta x^{(0)} / \lambda} \le \Delta \xi^{(k)} \le { \lambda \Delta x^{(0)}}$$
for an appropriate $\lambda$ (usually $\lambda=2$).

Further, by using the relation between $\psi(t)$ and $\varphi(t)$,
it is possible to incorporate the structure (\ref{112}) of the
s.-s.s. in the adaptation procedure.
In the case of equation (\ref{07}) we have
\be\label{2322}
\xi= x {\Gamma(t)}^m, \quad \Delta \xi=
\Delta x {\Gamma (t)}^m, \quad  m={{(\beta-\sigma-1)}/{2}}, \quad \Gamma (t)= \frac {\max_x u(t,x)} {\max_x u_0(x)}.
\ee
On the basis of the relations (\ref{2322}) the following strategy is accepted.

In the case of a single point blow-up, $m>0$, we choose the step  $\Delta x^{(k)}$
so that the step $ \Delta \xi^{(k)}$ be bounded from above:
\be\label{2323}
  \Delta \xi^{(k)} = \Delta x^{(k)} {\Gamma (t)}^m \leq { \lambda \Delta x^{(0)}}.
\ee
When condition (\ref{2323}) is violated, the following procedure is carried out:
 every element in the region, in which the solution is not established with a given accuracy $\delta_u$
(usually $\delta_u = 10^{-7}$), is divided into two equal elements and the values of the solution in the
new mesh points are found by interpolating from the old values;
the elements, in which the solution is established with a given accuracy $\delta_u$, are neglected.

In the case of a total blow-up, $m<0$, we choose the step  $\Delta x^{(k)}$
so that the step $ \Delta \xi^{(k)}$ be bounded from below:
\be\label{2324}
  \Delta \xi^{(k)} = \Delta x^{(k)} {\Gamma (t)}^m \geq { \lambda \Delta x^{(0)}}.
\ee
When condition (\ref{2324}) is violated, the lengths of the elements are doubled, and so is the interval in $x$:
$X(t_{k+1})=2X(t_k)$, thus the number of mesh points remains constant.

This adaptation procedure makes  it possible to compute efficiently the single point blow-up as well as
the total blow-up solutions up to amplitudes $10^6-10^{12},$ depending on the medium parameters.
It ensures the authenticity of the results of investigation of
the structural stability and the metastability of the self-similar solutions.  Let us note,
this approach does not require   an auxiliary differential problem for the mesh to be solved, unlike
the moving mesh methods  \cite{1100}. The idea to use the invariant properties
of the differential equations and their solutions \cite{114} and to incorporate   the structural properties
(e.g. geometry, different kind of symmetries,
the conservation laws) of the continuous problems in the numerical method, lies at the basis of
an important direction of the computational mathematics - geometric integration, to which
many works and monographs are devoted - see \cite{70}, \cite{68}, \cite{69}
and the references therein.

The numerous computational experiments carried out with the exact self-similar initial data (\ref{019}) for
$N=1, \beta=\sigma+1,$ as well as with the computed self-similar initial data, show a good
blow-up time restoration (set into the self-similar problem) in the process of solving the reaction-diffusion problem.
The preservation of the self-similarity and the restoration of the blow-up time
demonstrate the high quality  of the numerical methods for solving both the self-similar and
the nonstationary problems.

\section{Results and achievements}
\label{sec:4}
The developed numerical technique was used to analyze and solve a number of open problems. Below we present briefly
some of them.

\medskip
\noindent
{\bf 4.1 The transition $LS$- to $S$-regime in the radially symmetric case}

\noindent
The investigation of the limit case $\beta \to \sigma+1+0$ resolved the following paradox for
  $N>1$: there exists  one simple-structure s.-s. function in  $S$-regime $(\beta =\sigma+1),$
whereas in  $LS$-regime for
$\beta\to \sigma+1+0$ their number tends to infinity according to formula
(\ref{016}). The detailed numerical experiment in
\cite{D11}, \cite{D09} yielded the following results.
First, it was shown, that the structure of the s.-s.f. for $N>1$
and  $\beta\sim\sigma+1, \ \beta>\sigma+1$ is substantially different from the one for
$N=1.$ Second, for $N=1$ the transition $\beta\to
\sigma+1+0$ is ``continuous'' -- the self-similar function $\theta_{s,k}(\xi)$ for the $LS$-regime
tends to a s.-s.f. of the $S$-regime, consisting of  $k$ elementary solutions.
For $N>1$ the transition behaves very differently.
\begin{figure}
\begin{center}
\begin{tabular}{cc}
\includegraphics[width=0.4\textwidth]{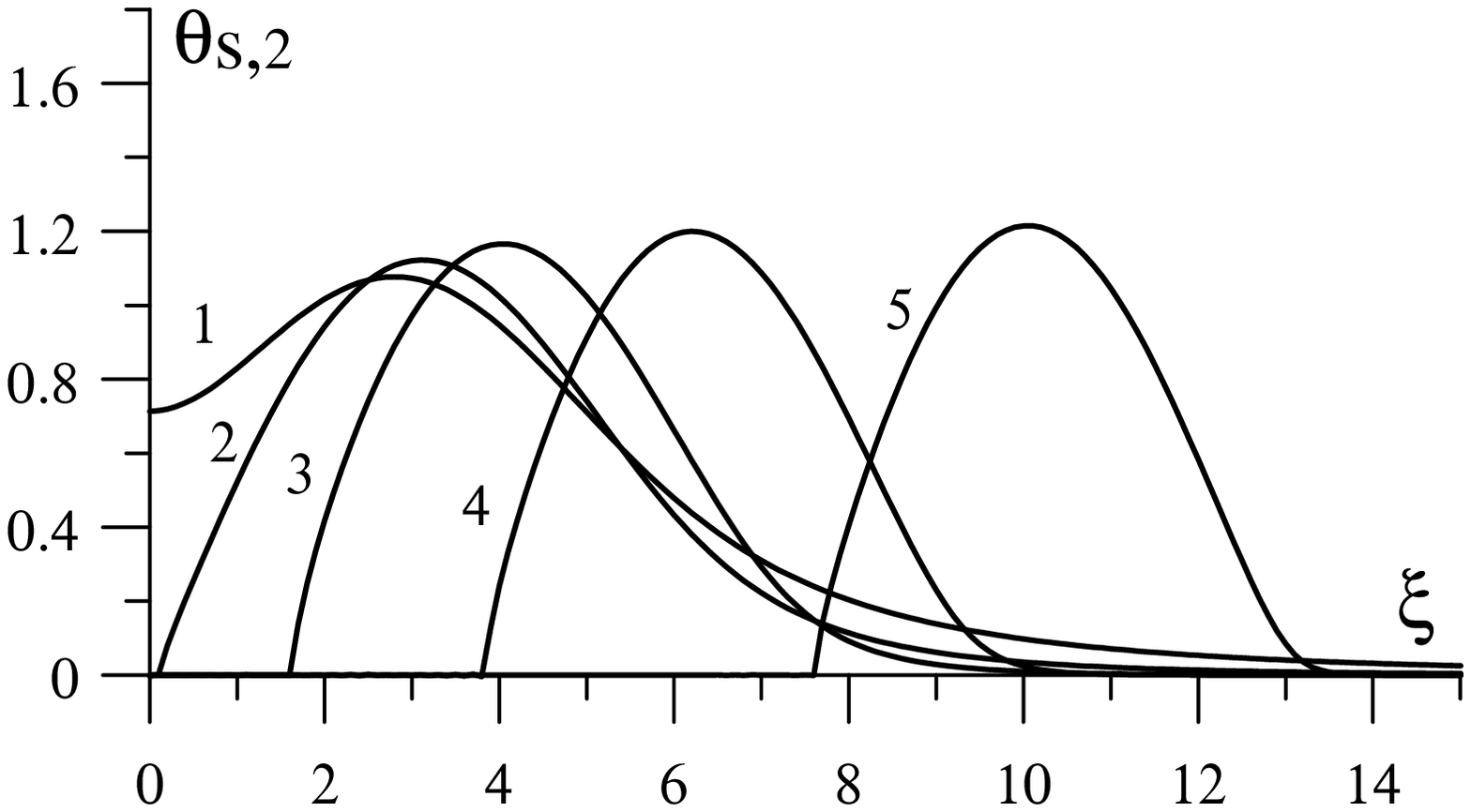} &
\includegraphics[width=0.4\textwidth]{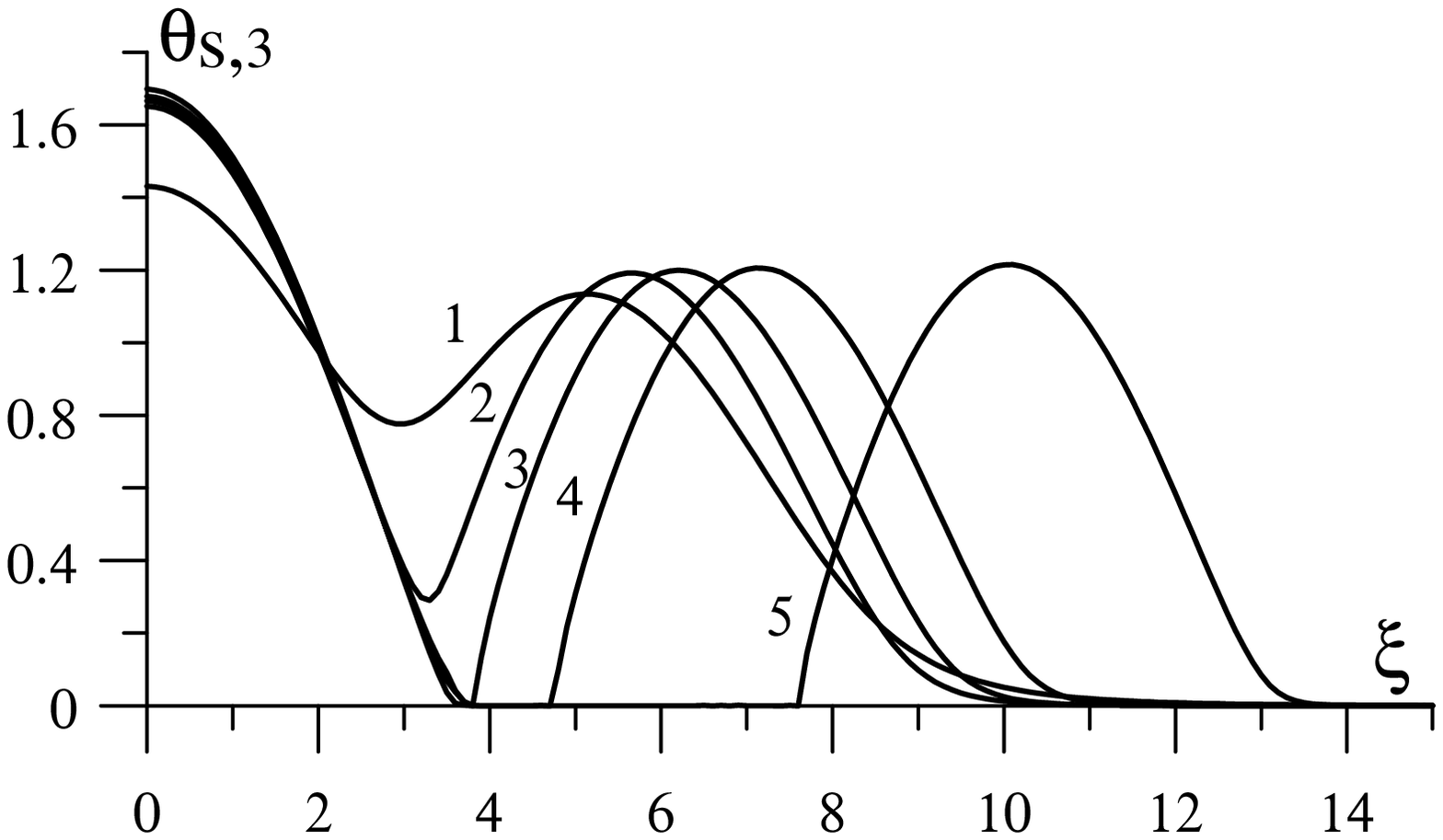}
\end{tabular}
\parbox[t]{0.9\textwidth}
{
\caption{\small {Graphs of the  s.-s.f. $\theta_{s,2}$ for $N=3, \ \sigma=2, \
\beta = \{3.6 (1); \ 3.38 (2); \  3.2 (3); \ 3.08 (4); \ 3.03(5)\}$ (left)
and $\theta_{s,3}$ for $N=3, \ \sigma=2, \
\beta= \{3.2 (1); \  3.1 (2); \ 3.08 (3); \ 3.06 (4); \ 3.03 (5)\}$ (right). }}}
\end{center}
\end{figure}

For $\sigma$ fixed and $\beta \to \sigma+1+0$
the central minimum of the "even"  s.-s.f. $\theta^{(N)}_{2j},
\ j=1,2,\ldots ,$ decreases and surprisingly becomes
 zero  for some $\beta=\beta_j^*(\sigma, N)$ (Fig. 5, left). For $\beta<\beta_j^*(\sigma, N)$ all
 of the s.-s.f. have zero region around the center of symmetry, and the radius of
 this region tends to infinity for $\beta \to \sigma+1+0$.
 All of the maxima of
$\theta^{(N)}_{s,2j}(\xi)$ tend to the maximum of the s.-s.f. of the   $S$-regime
for the corresponding  $\sigma$ and for $N=1.$ Thus the s.-s.f.
$\theta^{(N)}_{s,2j}(\xi),$ ``going to infinity'' when
$\beta\to\sigma+1+0$, tends to a s.-s.f. of the  $S$-regime for $N=1$ and the same
 $\sigma$, consisting of  $j$ elementary solutions.

\medskip

 For fixed  $\sigma$ there exists such a value $\beta^{**}_j(\sigma, N)$, that for
$\sigma+1<\beta<\beta^{**}_j$ the "odd" s.-s.f. $\theta^{(N)}_{s,2j+1}(\xi),
j=1,2,\ldots$ split into two parts: a central one, tending to the s.-s.f. of the
$S$-regime for the same $N$,  and second one, coinciding with the s.-s.f.
$\theta^{(N)}_{s,2j}(\xi),$ ``going to infinity'' when $\beta
\to \sigma+1+0$ (Fig. 5, right).

According to the described ``scenario'', when
$\beta\to\sigma+1+0$ only the first s.-s.f. of the $LS$-regime remains,
and it tends to the unique s.-s.f. of the $S$-regime.

As a result of this investigation new-structure s.-s.f. were found -- s.-s.f. with
{\bf a left front}. The existence of such s.-s.functions was confirmed by an asymptotic analysis, i.e.,
the asymptotics in the neighborhood  of the left front-point was found analytically \cite{D25}.
This new type of solutions  initiated investigations of other  authors  \cite{74},  \cite{132}
by other methods  (the method of dynamical analogy, bifurcation analysis). Their investigations
confirmed our results.

\medskip

\noindent
{\bf 4.2 The asymptotic behavior of the blow-up solutions of problem
 (\ref{07}), (\ref{08}) beyond the critical Fujita exponent}

\noindent
 For  $\beta> \beta_f = \sigma+1+{2/N}$ the problem
 (\ref{07}), (\ref{08}) could have blow-up- or global solutions depending on the initial data.
 For the s.-s. blow-up solution (\ref{09}), (\ref{010}) it holds $u_s(t,r)\not\in
L_1(\mathbb{R}^N)$.  The qualitative theory of nonstationary averaging ``amplitude-semi-width''
predicts a self-similar behavior of the amplitude of the blow-up solutions  and a possible
non-self-similar behavior of the semi-width \cite{32}. A question was posed there:
what kind of invariant or approximate s.-s.s. describes the asymptotic stage  $(t\to T^-_0)$ of the blow-up process?

The detailed  numerical experiment carried out in  \cite{D19}, \cite{78} showed that the s.-s.s.
(\ref{09}), (\ref{010}), corresponding to $\theta_{s,1}(\xi)$, is structurally stable:
all of the numerical experiments with finite-support initial data
(\ref{08}), ensuring blow-up, yield solutions tending to the self-similar one on the asymptotic stage.

\medskip

\noindent
{\bf 4.3 Asymptotically  self-similar blow-up  beyond
some other critical exponents}

\noindent
 The numerical investigation of the blow-up processes
in the radially symmetric case for high space-dimensions $N$ was carried out in \cite{D22}, \cite{D23}.
The aim was to check some hypotheses \cite{50} about the solutions of the s.-s. problem (\ref{011}), (\ref{012})
and about the asymptotic stability of the corresponding s.-s. solutions
for parameters beyond the  following critical exponents:

\noindent
$ \beta_s = (\sigma + 1)(N+2)/(N-2), \  N \geq 3$ \ (Sobolev's exponent);

\noindent
$\beta_u = (\sigma + 1)(1 + 4/(N-4-2 \sqrt {N-1})), \ N \geq 11; $

\noindent
$\beta_p = 1+ {3(\sigma + 1) + (\sigma^2 (N-10)^2 + 2 \sigma (5
\sigma + 1)(N-10) +9(\sigma + 1)^2)^{1/2}}/(N - 10),$   $N \geq 11.$

Self-similar functions, monotone in space, were constructed numerically for all of these cases,
thus confirming the hypotheses of their existence (not proved for $\beta>\beta_u$).
It was also shown, that the corresponding s.-s.s. are structurally stable,
thus confirming another   hypothesis  of \cite{50}.
Due to the strong singularity at the origin, the nonsymmetric Galerkin method and the special refinement
of the finite element mesh \cite{D23} were crucial for the success of these investigations.

\medskip

\noindent
{\bf 4.4 Two-component nonlinear medium}

\noindent
The methods, developed for the radially symmetric problems (\ref{07}),
(\ref{08}) and  (\ref{011}), (\ref{012}) were generalized in \cite{D20}, \cite{136}, \cite{78}
for the case of two-component nonlinear medium,
described by the system:
\be \left|
\begin{array}{l}
\displaystyle u_{1t}={1\over x^{N-1}}(x^{N-1}u_1^{\sigma_1}
u_{1x})_x+u_1^{\beta_1}u_2^{\gamma_2}, \quad x\in \mathbb{R}^1_+, \ N=1,2,3,\\[10pt]
\displaystyle u_{2t}={1\over x^{N-1}}(x^{N-1}u_2^{\sigma_2}
u_{2x})_x+u_1^{\gamma_1}u_2^{\beta_2}, \quad \sigma_i> 0, \ \beta_i>1, \
\gamma_i\ge 0, i=1,2.\end{array}
\right.
\label{025}
\ee
This system admits blow-up s.-s.s. of the form
\be
\begin{array}{ll}
\displaystyle u_{1s}=(1-{t /T_0})^{m_1}\theta_{1s}(\xi), \ \
&\xi=x/(1-{t / T_0})^n,\\[10pt]
\displaystyle u_{2s}=(1-{t / T_0})^{m_2}\theta_{2s}(\xi), \ \
&m_i<0, \ i=1,2,\end{array}
\label{026}
\ee
where
$$m_i={\alpha_i\over p},\ \ \alpha_i=\gamma_i+1-\beta_i,\ i=1,2, \ \ \
p=(\beta_1-1)(\beta_2-1)-\gamma_1\gamma_2,$$
$$n={m_1\sigma_1+1\over 2} = {m_2\sigma_2+1\over 2},\ \ \
\sigma_1(\gamma_2+1-\beta_2)=\sigma_2(\gamma_1+1-\beta_1).$$
The s.-s.f. satisfy the system of nonlinear ODE:
\be
\left|
\begin{array}{l}
L_1(\theta_{1s}, \theta_{2s})\equiv -{1\over
\xi^{N-1}}(\xi^{N-1}\theta_{1s}^{\sigma_1}\theta'_{1s})'+n\xi\theta'_{1s} - m_1\theta_{1s}-
\theta_{1s}^{\beta_1}\theta_{2s}^{\gamma_2}=0,\\[10pt]
L_2(\theta_{1s}, \theta_{2s})\equiv -{1\over
\xi^{N-1}}(\xi^{N-1}\theta_{2s}^{\sigma_2}\theta'_{2s})'+n\xi\theta'_{2s} - m_2\theta_{2s}
-\theta_{1s}^{\gamma_1}\theta_{2s}^{\beta_2}=0\end{array}
\right.\label{027}
\ee
and the boundary conditions
\be \lim_{\xi \to 0} \xi^{N-1}\theta_{is}^{\sigma_1}\theta'_{is}=0,\ \ \lim_{\xi\to
\infty}\theta_{is}=0,\ \ i=1,2.\label{028}
\ee
Superconvergence of the FEM (of order  $O(h^4)$)
for solving the s.-s. problem (\ref{027})--(\ref{028}) by means of quadratic elements
and optimal-order convergence ($O(h^2)$) by means of linear elements were achieved.
The structural stability of the s.-s.s. (\ref{026}) for parameters
 $\sigma_i,\ \beta_i, \ \gamma_i, \ i=1,2,$ corresponding to the $LS$-regime ($n>0$)
 was analyzed in \cite{136}, \cite{78}, \cite{D20}.
It was shown that only the s.-s.s. of systems (\ref{025}) with strong feedback ($p<0$),
corresponding to the s.-s.f. with two simple-structure components, were structurally stable.
All the other s.-s.s. were metastable -- self-similarity was preserved to times not less than
 $99.3\% \tilde T_0$.

 The proposed computational technique can be applied to investigating the self-organization processes
 in wide classes of nonlinear dissipative media described by nonlinear reaction-diffusion systems.

\medskip

{\bf 4.5 Directed heat diffusion in a nonlinear anisotropic medium}

\noindent
Historically the first Bulgarian contribution to the topic under consideration was
the numerical realization of the  self-similar solutions, describing
directed heat diffusion and burning of a two-dimensional nonlinear anisotropic medium.
It was shown in  \cite{27} that the model of heat structures in the anisotropic case
\be
u_t=(u^{\sigma_1} u_{x_1})_{x_1}+(u^{\sigma_2}u_{x_2})_{x_2} +u^{\beta},\ x=(x_1,
x_2)\in \mathbb{R}^2, \ \sigma_1>0,\ \ \sigma_2>0, \ \ \beta>1 \label{039}
\ee
admits  invariant solutions of the kind
$$u_s(t,x_1,x_2)=\left(1-{t\over T_0}\right)^{-{1\over \beta-1}} \theta_s(\xi), \
\xi=(\xi_1,\xi_2)\in \mathbb{R}^2,$$
$$\xi_i=x_i / \left(1-{t\over T_0}\right)^{m_i\over{\beta-1}},\ \ m_i={\beta-\sigma_i-1\over 2},
\ i=1,2.$$
The self-similar function $\theta_s(\xi_1,\xi_2)$ satisfies the nonlinear elliptic problem
\be
L(\theta_s)\equiv \sum^2_{i=1}\left(-{\partial\over \partial
\xi_i}\left(\theta_s^{\sigma_i}{\partial \theta_s\over \partial\xi_i}\right)
+{\beta-\sigma_i-1\over 2} \xi_i{\partial\theta_s\over
\partial\xi_i}\right)+\theta_s-\theta_s^{\beta}=0,\label{041}
\ee
\be
{\partial\theta_s\over \partial\xi_i}\Biggl.\Biggr|_{\xi_i=0}=0,\ i=1,2;\ \
\theta_s(\xi)\to 0,\ \ |\xi|\to \infty. \label{042}
\ee
The Cauchy problem for equation
 (\ref{039}) was investigated in the works \cite{D01},  \cite{D03} for different parameters $\sigma_1,
\sigma_2$ and $\beta$. Depending on the parameters,
 different mixed regimes: $S-HS, \ HS-LS, \ S-LS$ of heat transfer and burning
were implemented  numerically.
The evolution in time of one and the same initial perturbation is shown
for the cases of the 2D radially symmetric $S$-regime (Fig. 6), the mixed $HS-S$-regime (Fig. 7)
and  the mixed $HS-LS$-regime (Fig. 8).
\begin{figure}
\begin{center}
\includegraphics[width=0.8\textwidth]{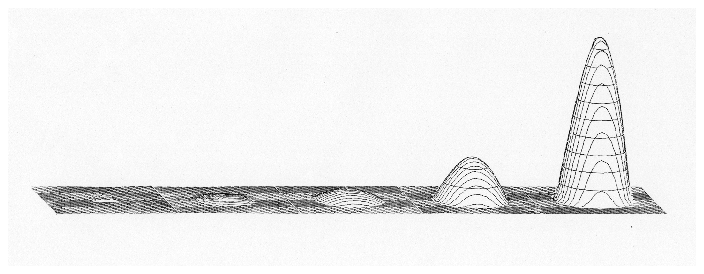}  \\
\end{center}
\caption{\small{$S$-regime: $\sigma_1=2,\ \sigma_2=2,\ \beta=3. $}}
\end{figure}
\begin{figure}
\begin{center}
\includegraphics[width=0.8\textwidth]{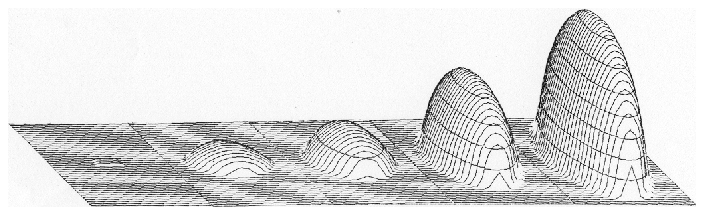}  \\
\end{center}
\caption{\small{$HS-S$-regime: $\sigma_1=3,\ \sigma_2=2,\ \beta=3. $}}
\end{figure}
\begin{figure}
\begin{center}
\includegraphics[width=0.8\textwidth]{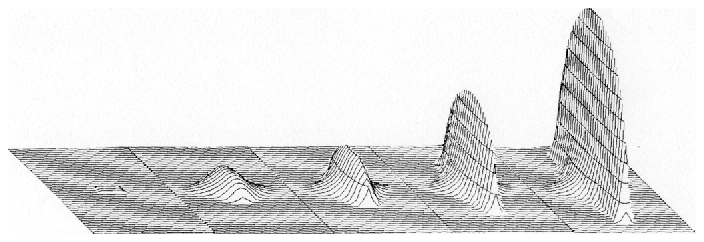}  \\
\end{center}
\caption{\small{$HS-LS$-regime: $\sigma_1=3,\ \sigma_2=1,\ \beta=3. $}}
\end{figure}

To solve the Cauchy problem for equation (\ref{039}) a modification of the TERMO Package of Applied Programs \cite{104},
designed initially for solving isotropic problems with piecewise constant coefficients, was done. TERMO had been worked out
by an IMI-BAS team, after the idea of
Raytcho Lazarov and under his guidance, a merit worth mentioning here.

 The s.-s.functions for the corresponding mixed regimes were found in \cite{D01} by self-similar processing
 of the solution of the Cauchy problem for equation (\ref{039}).
Later, in \cite{D06}, they were found as solutions of the self-similar problem (\ref{041}), (\ref{042}).
 The self-similar functions of complex symmetry for the isotropic case
in Cartesian  coordinates (denoted in \cite{57} as $Ei/j$) were found  as a special case.

 Later on, in \cite{D14}, \cite{136}, the numerical methods were modified for the isotropic 2D self-similar problem
   in polar coordinates to construct numerically  another class of self-similar functions of complex
   symmetry  (denoted in \cite{57} as  $EjMm$)  in  $LS$-regime  and to investigate their structural stability.

Graphical representations of the evolution of the anisotropic  invariant solutions, as well as
the s.-s. functions $EjMm$ for some different values of $j, m, \sigma, \beta$, are included in the Handbook \cite{85}.
We show some of the s.-s.f. $EjMm$ in Fig. 9.
\begin{figure}
\begin{center}
\begin{tabular}{ccc}
\includegraphics[width=0.25\textwidth,bb=67 606 220 760]{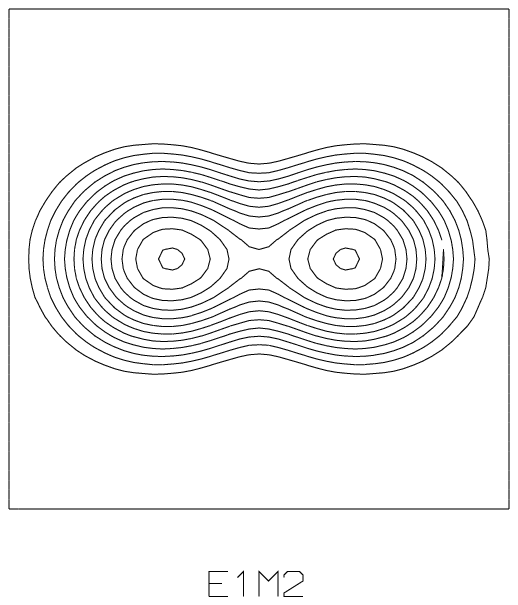} &
\includegraphics[width=0.25\textwidth,bb=320 606 473 760]{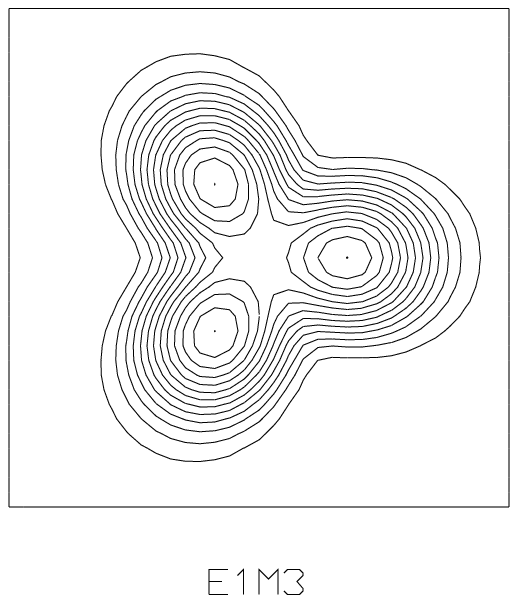} &
\includegraphics[width=0.25\textwidth,bb=67 373 220 526]{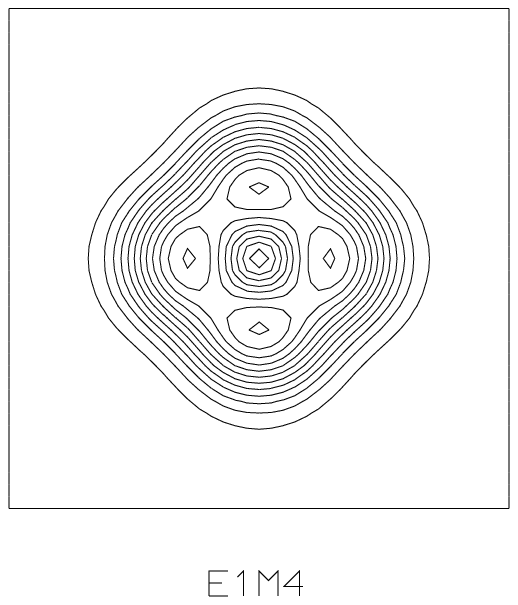}
\end{tabular}
\end{center}
\begin{center}
\begin{tabular}{ccc}
\includegraphics[width=0.25\textwidth,bb=320 373 473 526]{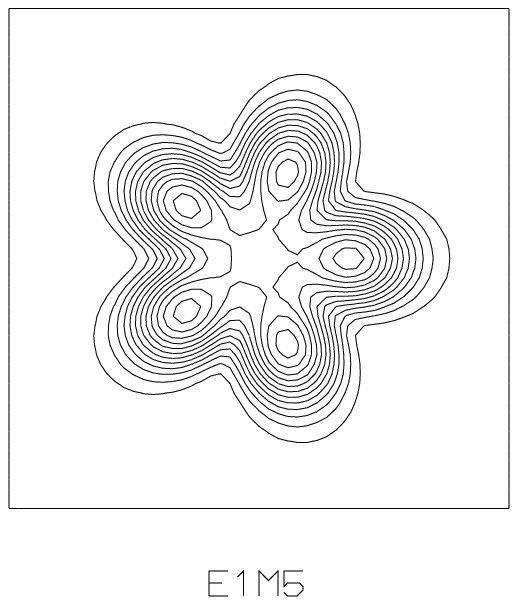} &
\includegraphics[width=0.25\textwidth,bb=67 140 220 293]{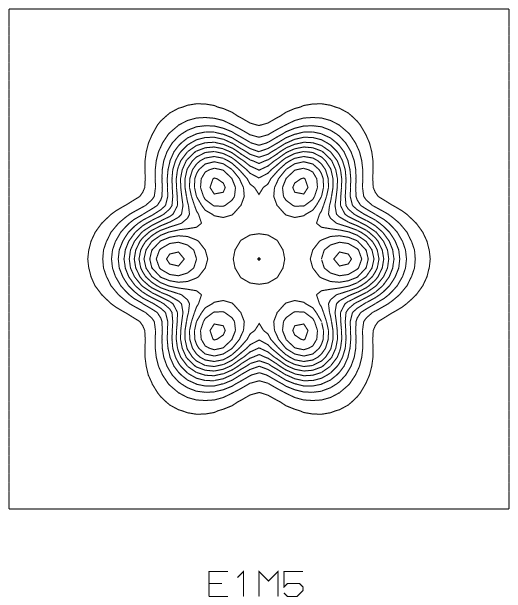} &
\includegraphics[width=0.25\textwidth,bb=320 140 473 293]{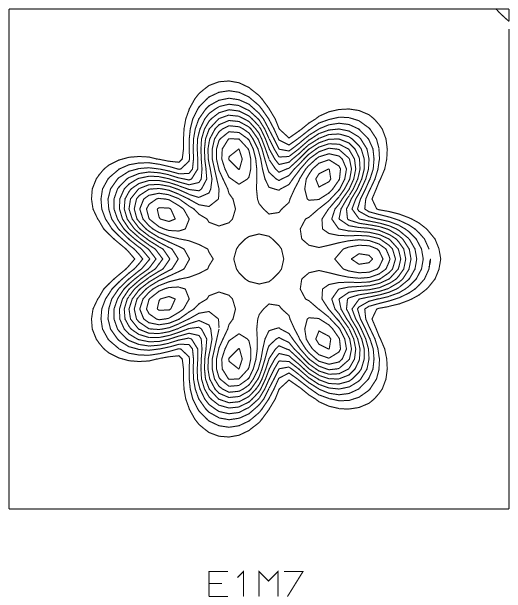}
\end{tabular}
\end{center}
\caption{\small{Self-similar functions EjMm, $j=1, \ m=2,3,4,5,6,7, \ $
$\sigma=2, \ \beta=3.25$.\label{fig111}}}
\end{figure}

\medskip

\noindent
{\bf 4.6 Spiral waves in $HS$-regime}

\noindent
The numerical realization of the invariant solutions,
describing  ``spiral'' propagation  of the nonhomogeneities in two-dimensional isotropic medium appears to be one
 of the most interesting contributions of ours.
The mathematical model in polar coordinates reads:
\be
u_t={1\over r}(r u^{\sigma} u_r)_r +{1\over
r^2}(u^{\sigma}u_{\varphi})_{\varphi}+u^{\beta}, \ \sigma>0,\
\beta>1.\label{030}
\ee
It admits s.-s.s. of the kind \cite{D03}, \cite{33}:
\be
u_s(t,r,\varphi)=\left(1-{t\over T_0}\right)^{-{1\over
\beta-1}}\theta_s(\xi,\phi),\label{031}
\ee
\be \xi = r / \left(1-{t\over T_0}\right)^{m\over{\beta-1}}, \ \ \phi =\varphi+{c_0\over
\beta-1}\ln \left(1-{t\over T_0}\right), \ \
m={\beta-\sigma-1\over 2}.\label{032}
\ee
The self-similar function $\theta_s(\xi,\phi)$ satisfies  the nonlinear elliptic equation
\be
\renewcommand{\arraystretch}{2.3}
\begin{array}{ll}
\ds L(\theta_s)\equiv &- \ds {1\over \xi} \ds{\partial\over \partial\xi}\left(\xi
\theta_s^{\sigma} {\partial\theta_s\over
\partial\xi}\right)-{1\over \xi^2}{\partial\over
\partial\phi} \left(\theta^{\sigma}_s{\partial\theta_s\over
\partial\phi}\right)+{\beta-\sigma-1\over 2} \xi {\partial\theta_s\over
\partial\xi}\\
\ds & - \ds c_0{\partial\theta_s\over \partial\phi} +\theta_s-\theta_s^{\beta}=0,\ \ \
T_0={1\over \beta-1}. \end{array} \label{033}
\renewcommand{\arraystretch}{1}
\ee
Here $c_0\not = 0$ is the parameter of the family of  solutions. From (\ref{032})
it follows
$$\xi e^{s\phi}=r e^{s\varphi}={\rm const}, \ \ s={\beta-\sigma-1\over
2c_0}.$$
This means, that the trajectories of the nonhomogeneities in the medium (say local maxima)
are logarithmic spirals for  $\beta\not = \sigma+1$ or circles for
 $\beta=\sigma+1.$ The direction of movement for fixed $c_0$,
for example $c_0>0,$ depends on the relation between $\sigma$ and $\beta:$  for
$\beta > \sigma+1$ -- towards the center (twisting spirals), for $\beta<\sigma+1$ --
from the center (untwisting spirals).

 The problem for the numerical realization of the spiral s.-s.s. (\ref{031}), (\ref{032}) was posed in 1984, when
the possibility for their existence has been  established
 by the method  of invariant group analysis in the PhD Thesis of S.R. Svirshchevskii.
 As it was stated in \cite{82},  there were significant difficulties for finding such solutions.
 First, the linearization of the self-similar equation  (\ref{033}) was not expected to give
 the desired result, because it is not possible to separate the variables
 in the linearized  equation. Second, the asymptotics at infinity of the solutions of
 the self-similar equation were not known.

The first successful step was the appropriate (complex) separation of variables in the linearized  equation.
Using the assumption for small oscillations of the s.-s.f. $\theta_s(\xi,\phi)$ around
 $\theta_H^1\equiv 1$, i. e.,
$ \theta_s(\xi,\phi)=1+\alpha y(\xi,\phi),\ \ \alpha=\mbox{const}, \ \  |\alpha y|\ll 1$
and the idea of linearization around it,  the following
linear equation for $y(\xi,\phi)$ was found \cite{D16}:
$$
 -\frac{1}{\xi} \frac{\partial}{\partial\xi}\left(\xi
\frac{\partial y}{\partial\xi}\right) -\frac{1}{\xi^2}
\frac{\partial^2y}{\partial\phi^2}
+\frac{\beta-\sigma-1}{2}\xi\frac{\partial y}{\partial\xi}
-c_0\frac{\partial y}{\partial\phi} +(1-\beta)y=0. $$
Seeking for particular solutions $ Y_k(\xi, \phi)=R_k(\xi)e^{\rm{i}\phi},$   $k \in \mathbb{N},$
it was found: for $\beta=\sigma+1$
$$R_k(\xi)=J_k(z), \ \ \ z=(\sigma + c_0k\rm{i})^{1/2}\xi,$$
where $J_k(z)$ is the first kind Bessel function of order $k$,
and for $\beta\not=\sigma+1$
$$R_k(\xi)=\xi^k {_1F_1} (a,b;z), \quad
a=-{\beta-1+c_0 k\rm{i}\over \beta -\sigma -1} +{k\over 2}, \ \  b=1+k, \ \
z={\beta-\sigma-1\over 4}\xi^2,$$
where $_1 F_1(a,b,z)$ is the confluent hypergeometric function.
  The detailed analytical and numerical investigation \cite{D16} of the functions $y_k(\xi,\phi)=\Re(Y_k(\xi, \phi))$
 showed  that
their asymptotics at infinity are self-similar,  as well as that the functions
\be
\tilde{\theta}_{s,k}(\xi,\phi)=1+\alpha y_k(\xi,\phi),\ \ |\alpha y_k| \ll 1 \label{eq5}
\ee
are very close to the sought-after solutions $ \theta_s(\xi,\phi).$
Moreover, the amplitude of the linear approximations $y_k(\xi,\phi)$  tends to zero for
$\xi\to \infty$ in the case of $HS$-regime, and to infinity in the case of $LS$-regime. This gave
the idea to seek for s.-s.s. of the $HS$-regime,  tending to
the nontrivial constant solution  $\theta_s \equiv \theta_H,$  i.e., to generalize the notion of s.-s.
functions,  and consequently, the notion of the structures and waves, which arise and preserve themselves in the
absolutely cold medium. Although not exploited earlier, this change is reasonable and fully adequate to the real systems.
\begin{figure}
\begin{center}
\begin{tabular}{ccc}
\includegraphics[width=0.27\textwidth]{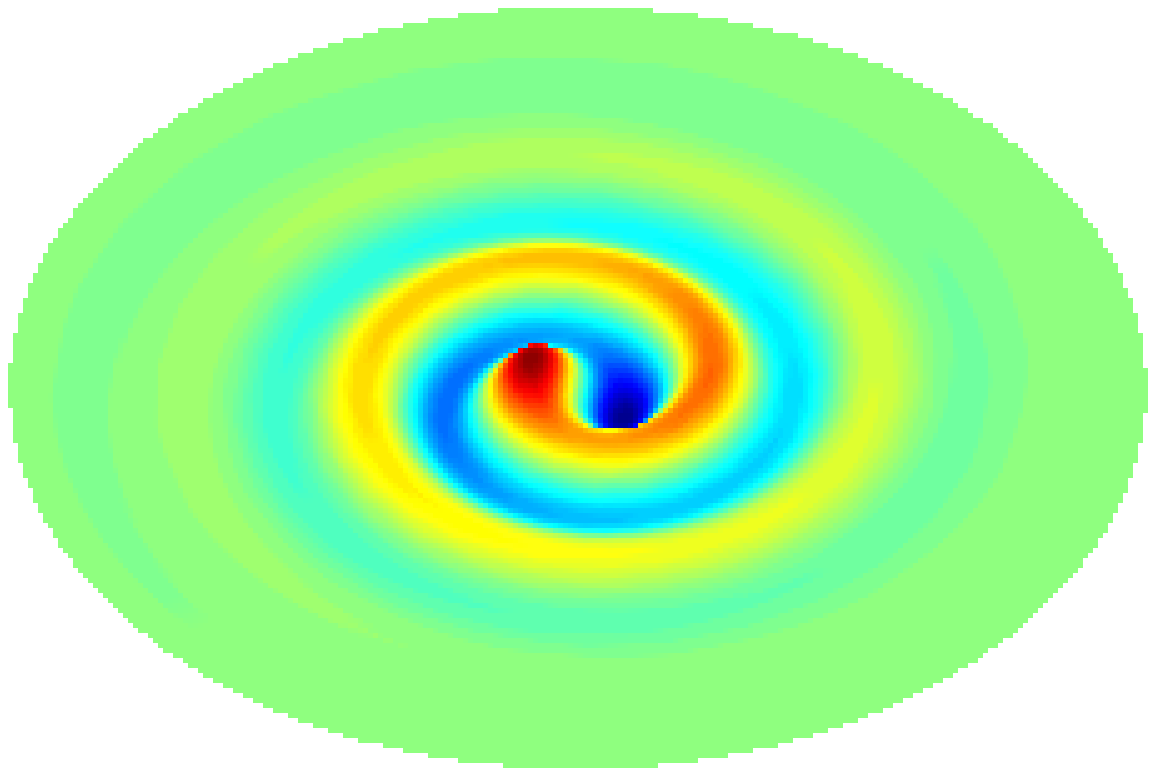} &
\includegraphics[width=0.27\textwidth]{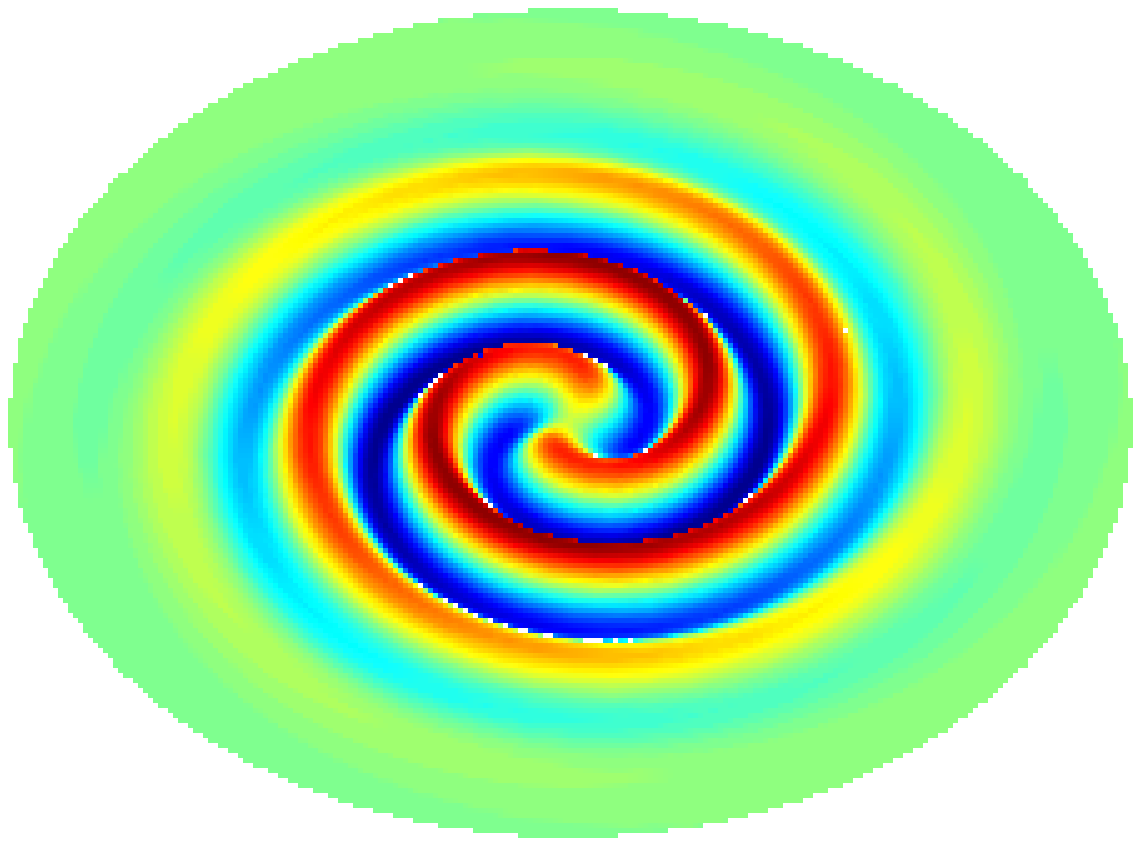} &
\includegraphics[width=0.27\textwidth]{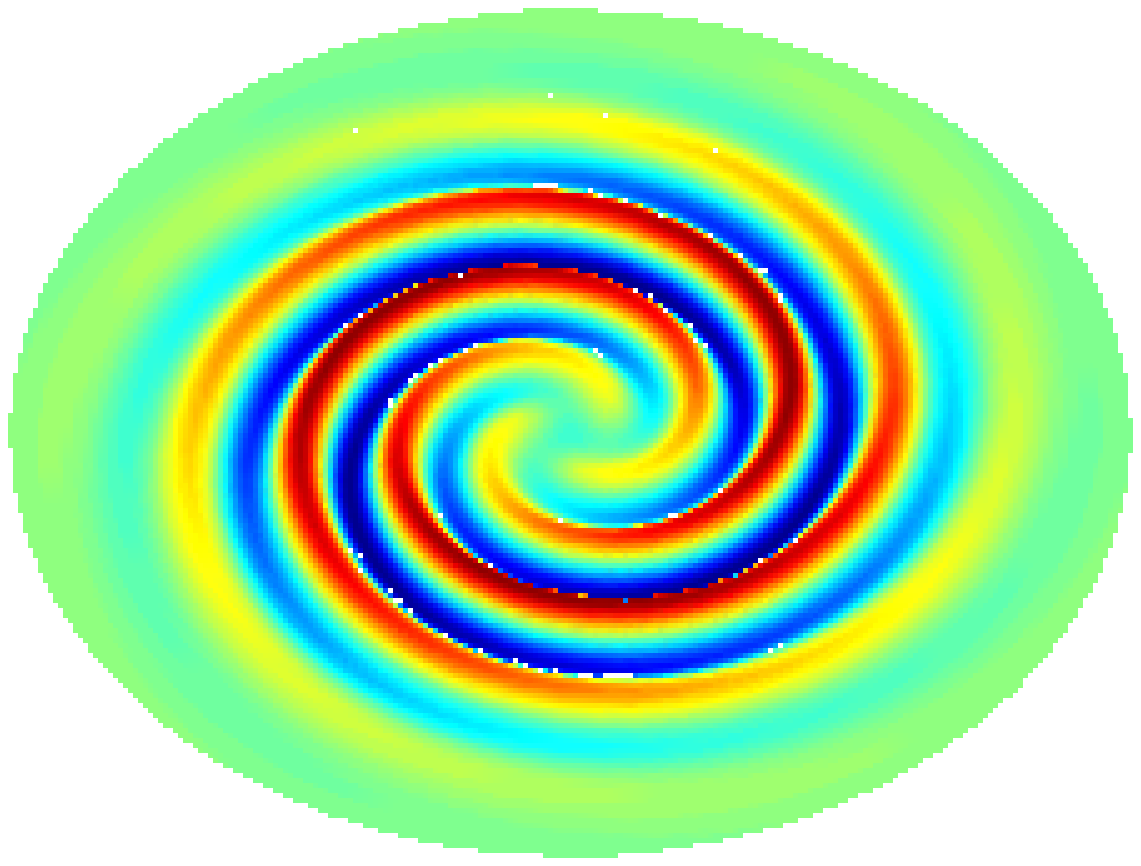}  \\
$k=1$ & $k=2$& $k=3$
\end{tabular}
\parbox[t]{0.85\textwidth}{
\caption{\small{One-armed spiral solution ($k=1$), two-armed spiral
solution ($k=2$), three-armed spiral solution ($k=3$),  $c_0=1$,
$\sigma=3,\beta=3.6$.\label{fig11}}}}
\end{center}
\end{figure}

All stated  above enabled us
to predict the asymptotics of the solutions of equation (\ref{033}),
to derive a boundary condition by using this asymptotics
and to close the s.-s. problem by the following boundary
and periodic conditions \cite{D24}, \cite{D18}, \cite{D21}:
\be
\renewcommand{\arraystretch}{2.3}
\begin{array}{l}
\ds\lim_{\xi\to 0}
\xi\theta_{s,k}^\sigma\frac{\partial\theta_{s,k}}{\partial\xi}=0,\
\phi\in\left[0,{2\pi \over k}\right],\\
  \ds \frac{\partial\theta_{s,k}}{\partial\xi}=\ds \frac{\theta_{s,k}-1}{\bar{m}
\xi} -\frac{\gamma k}{s \xi^{1 - \frac{1}{\bar{m}}}} \sin( k \phi + \frac{k}{s}\ln\xi
+ \mu), \ \xi=l\gg 1, \ \phi\in\left[0, {2\pi \over k}\right], \\
\ds \theta_{s,k}(\xi,0)=\theta_{s,k}\left(\xi,{2\pi \over k}\right), \ \ \
\ds\frac{\partial\theta_{s,k}}{\partial\phi}(\xi,0)=
\ds\frac{\partial\theta_{s,k}}{\partial\phi}\left(\xi,{2\pi \over k}\right), \ \ \
0\leq\xi\leq l,
\end{array}
\label{eq063q}
\ee
where $\bar{m}=m/(\beta-1)$, $\gamma$ and $\mu$ are constants, depending on $\sigma, \beta, c_0$.
The numerical solving of the problem  (\ref{033}), (\ref{eq063q}) for $c_0 \neq 0$ with initial approximations
(\ref{eq5})  gives  ``the spiral'' s.-s. functions of the $HS$-regime; some of them and their evolution in time
by solving equation  (\ref{030}) were investigated in   \cite{D18}, \cite{D21}, \cite{D24}.
The graphs of the s.-s.f. for $k=1$ (one-armed
spiral), $k=2$ (two-armed spiral), $k=3$ (three-armed spiral) are
shown  in Fig.~\ref{fig11}. The rest of the parameters are $\sigma=3$, $\beta=3.6$,
$c_0=1$.
\begin{figure}[t]
\begin{center}
\includegraphics[width=0.8\textwidth]{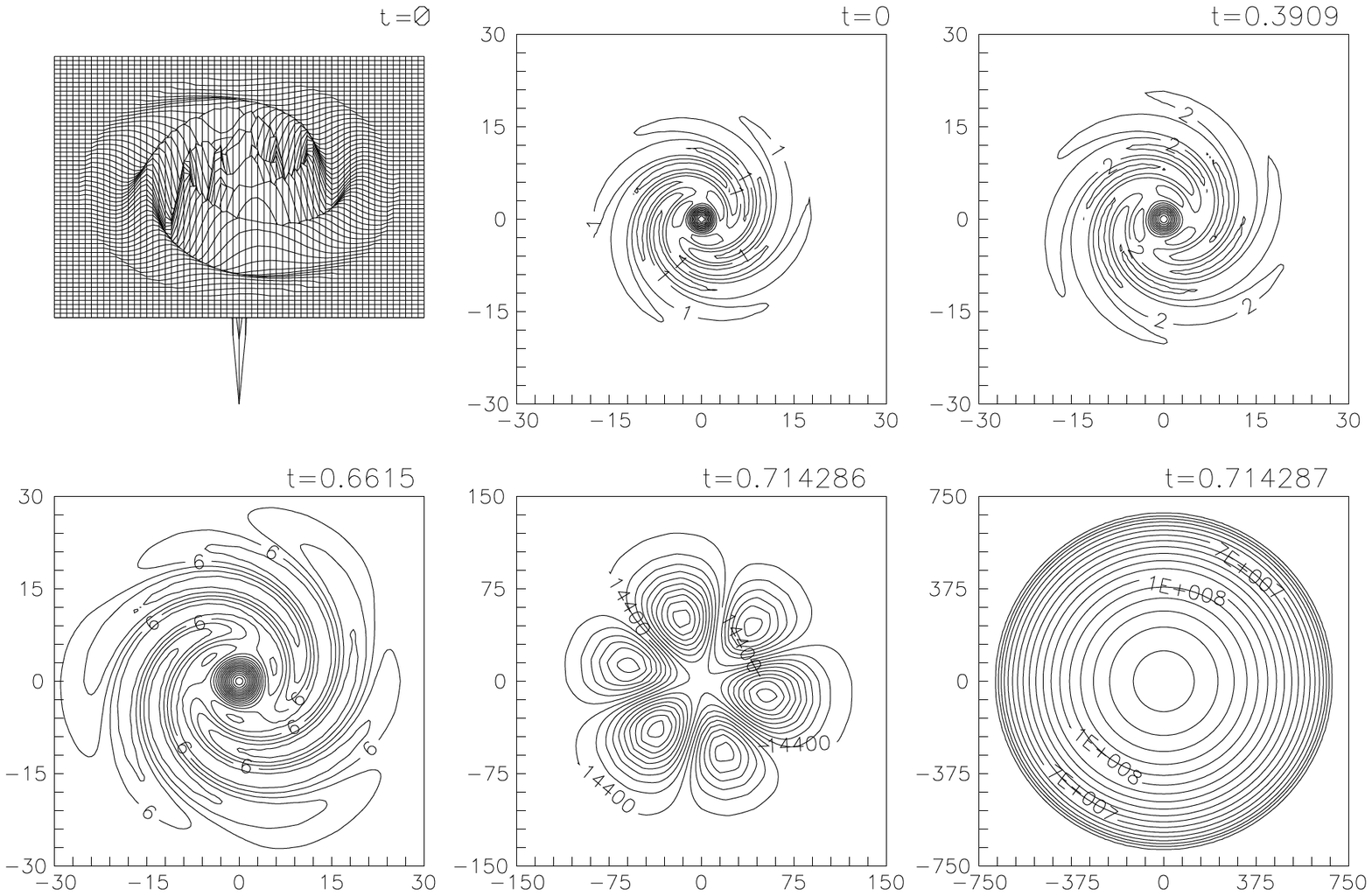}  \\
\parbox[t]{0.85\textwidth}{
\caption{\small{Evolution of three-armed spiral wave:
$\sigma=2,\ \beta=2.4,\ c_0=1,\ k=3,\ T_0= 0.(714285)$.\label{fig12}}}}
\end{center}
\end{figure}

The evolution of three-armed s.-s.f. for parameters $\sigma=2$, $\beta=2.4$, $k=3$,
$c_0=1$ is shown in Fig.~\ref{fig12}. The exact blow-up time is (\ref{014}) $T_0=1/(\beta-1)= 0.(714285).$
Similarly to  all complex-symmetry s.s.-s., the three-armed spiral one  is metastable -- at $T \to T_0$ it
degenerates into the simplest radially symmetric s.-s.s. for the same parameters $\sigma, \beta$.
The mesh adaptation in $r$-direction when solving equation (\ref{030}) is realized
again in consistency with the self-similar law,
keeping the same number of mesh points during the whole process of evolution.
Having to solve two nonlinear 2D problems (elliptic and reaction-diffusion),
going through a number of approximations, it is astonishing that the restoration (with accuracy $10^{-6}$)
of the exact blow-up time  is practically  perfect (Fig. 11, right most bottom).

The existence of spiral structures in
$LS$-regime is an open question by now. The ``ridges'' of the linear approximations of the s.-s.f. in the
 $LS$-regime  tend also to the self-similar ones for $\xi\to\infty$
 \cite{D16}, but their amplitudes tend to infinity and it is not clear by what asymptotics
 to sew the linear approximations.

\medskip

{\bf 4.7 Complex symmetry waves in $HS$-regime and $S$-regime}

\begin{figure}[t]
\begin{center}
\includegraphics[width=0.8\textwidth]{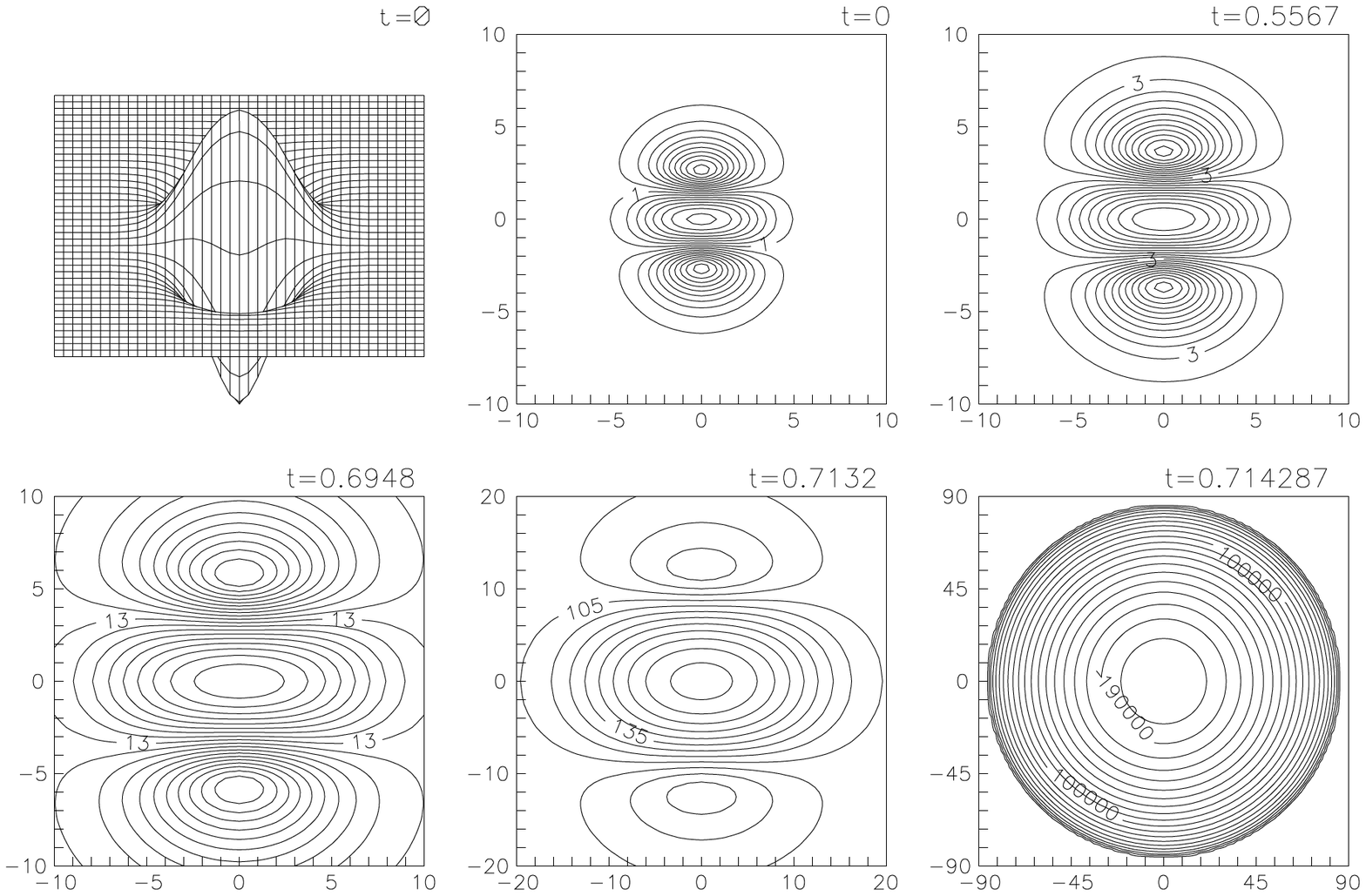}  \\
\parbox[t]{0.85\textwidth}{
\caption{\small{Evolution of a complex wave in $HS$-regime:
$\sigma=2,\ \beta=2.4,\ c_0=0,\ k=2,\ T_0= 0.(714285)$.\label{fig13}}}}
\end{center}
\vspace{-0.2in}
\end{figure}

\noindent
In the process of solving the problem for the spiral s.-s.s.  an idea arose to seek for
complex nonmonotone waves, tending to the nonzero homogeneous solution for $\xi \to \infty$  in
$HS$-regime and $c_0=0.$
Fig.~\ref{fig13} shows a complex-symmetry s.-s.f. and its evolution in time for the same parameters
$\sigma, \beta$, as in
Fig.~\ref{fig12}, and $k=2$, $c_0=0$. Note the same perfect restoration of the blow-up time.
The results about the spiral and the complex-symmetry s.-s.f. are included in the book   \cite{86}.

\begin{figure}[t]
\begin{center}
\begin{tabular}{ccc}
\includegraphics[width=0.27\textwidth]{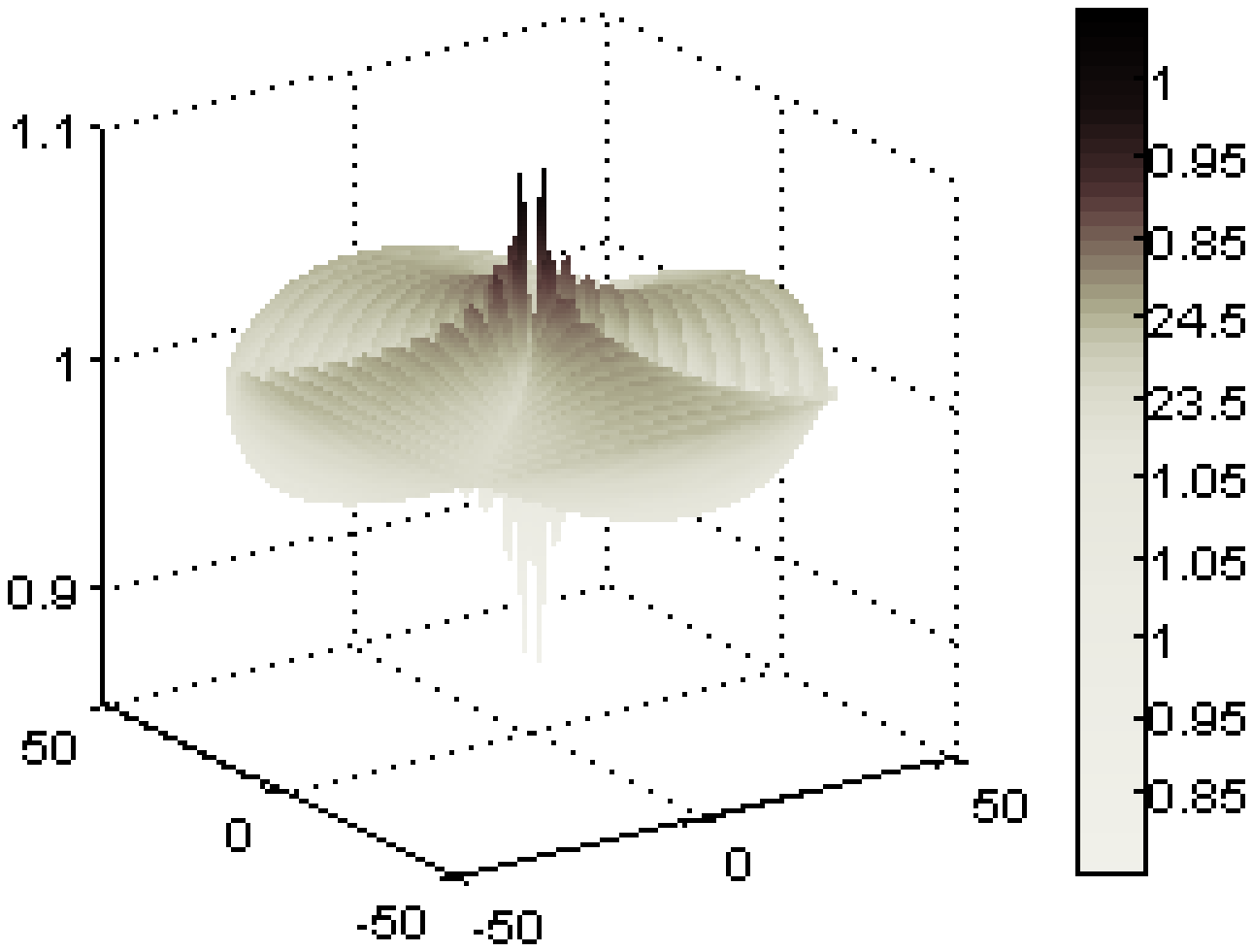} &
\includegraphics[width=0.27\textwidth]{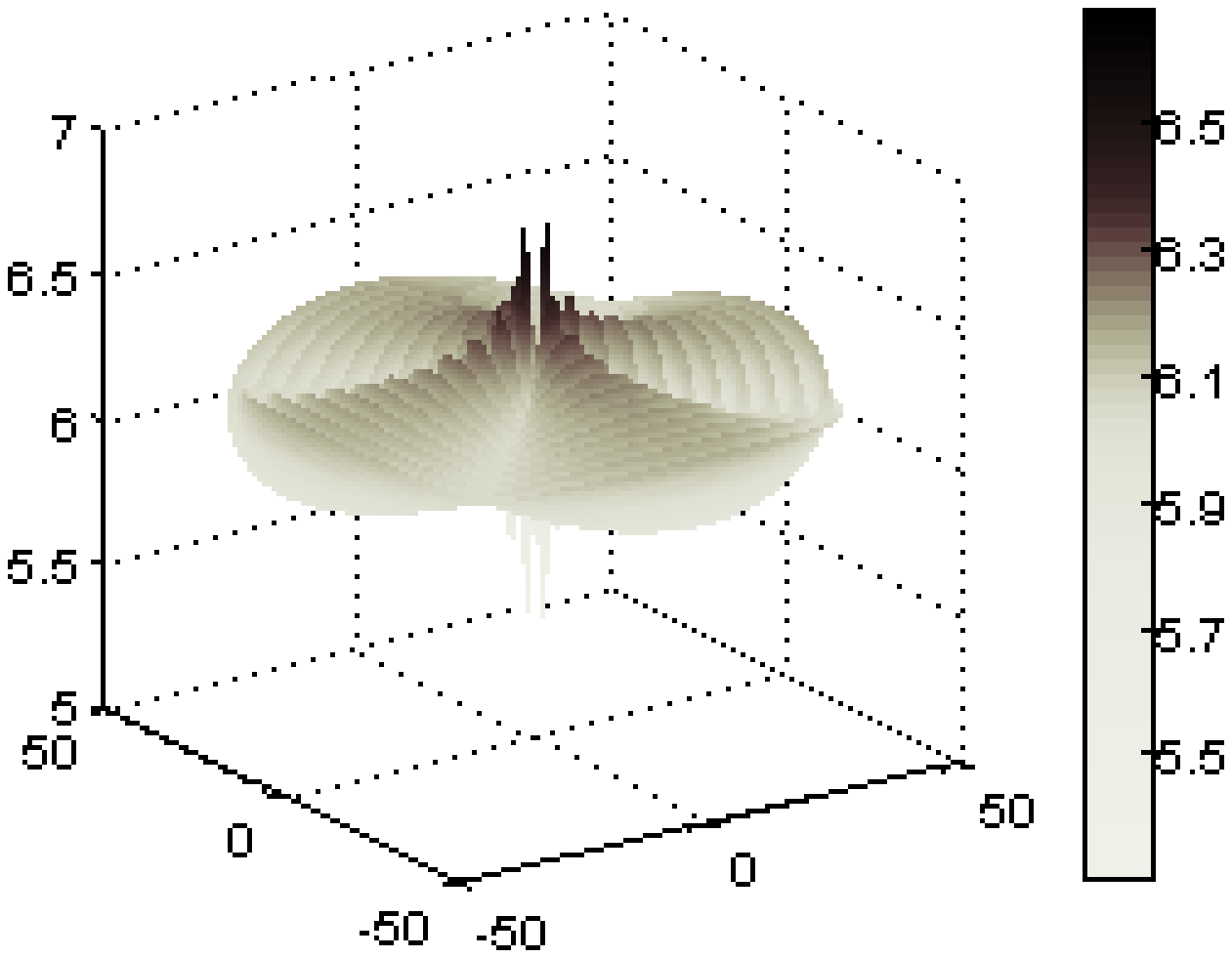} &
\includegraphics[width=0.27\textwidth]{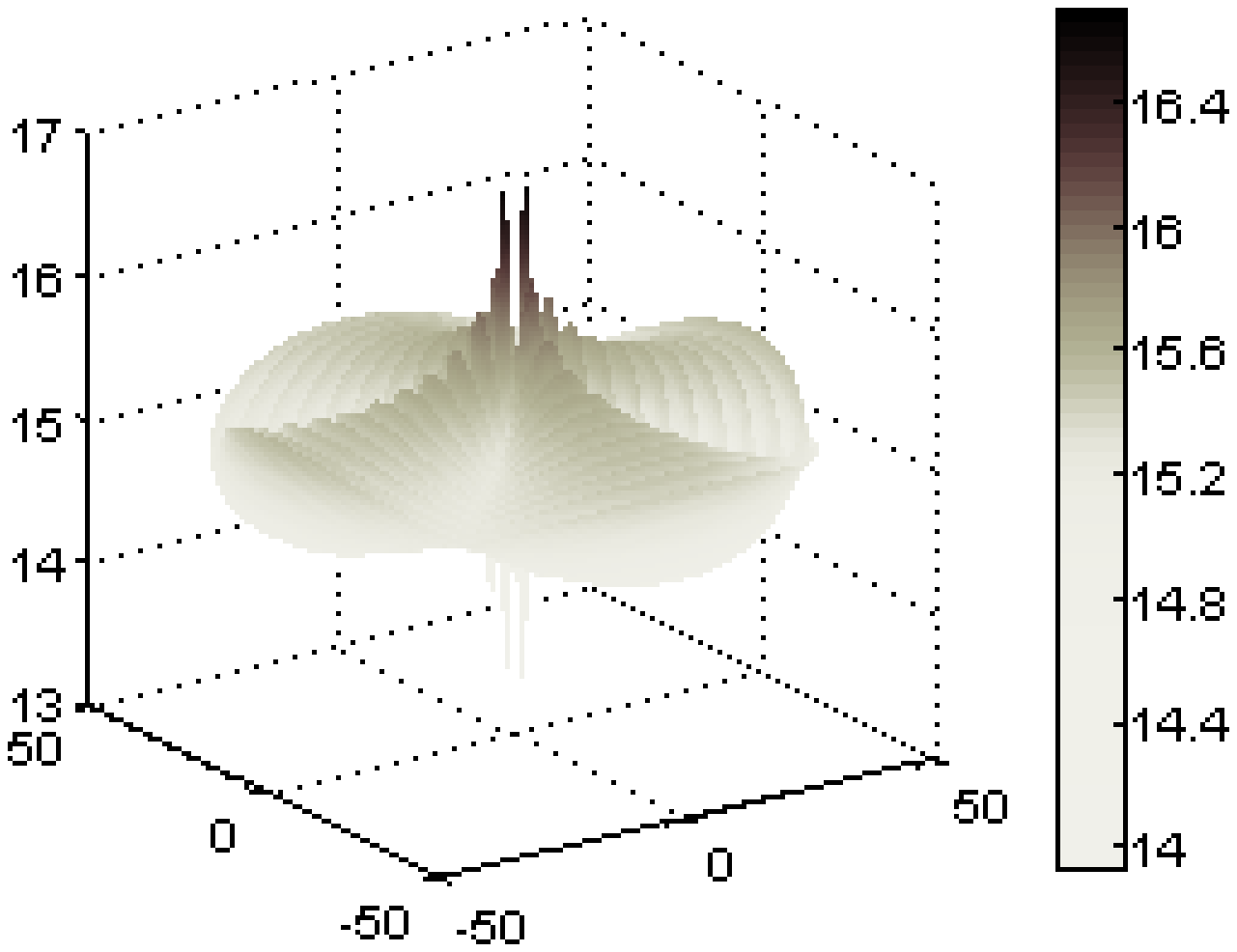} \\
$t=0$ & $t=0.486748$ & $t=0.497776$  \\[7mm]
\includegraphics[width=0.27\textwidth]{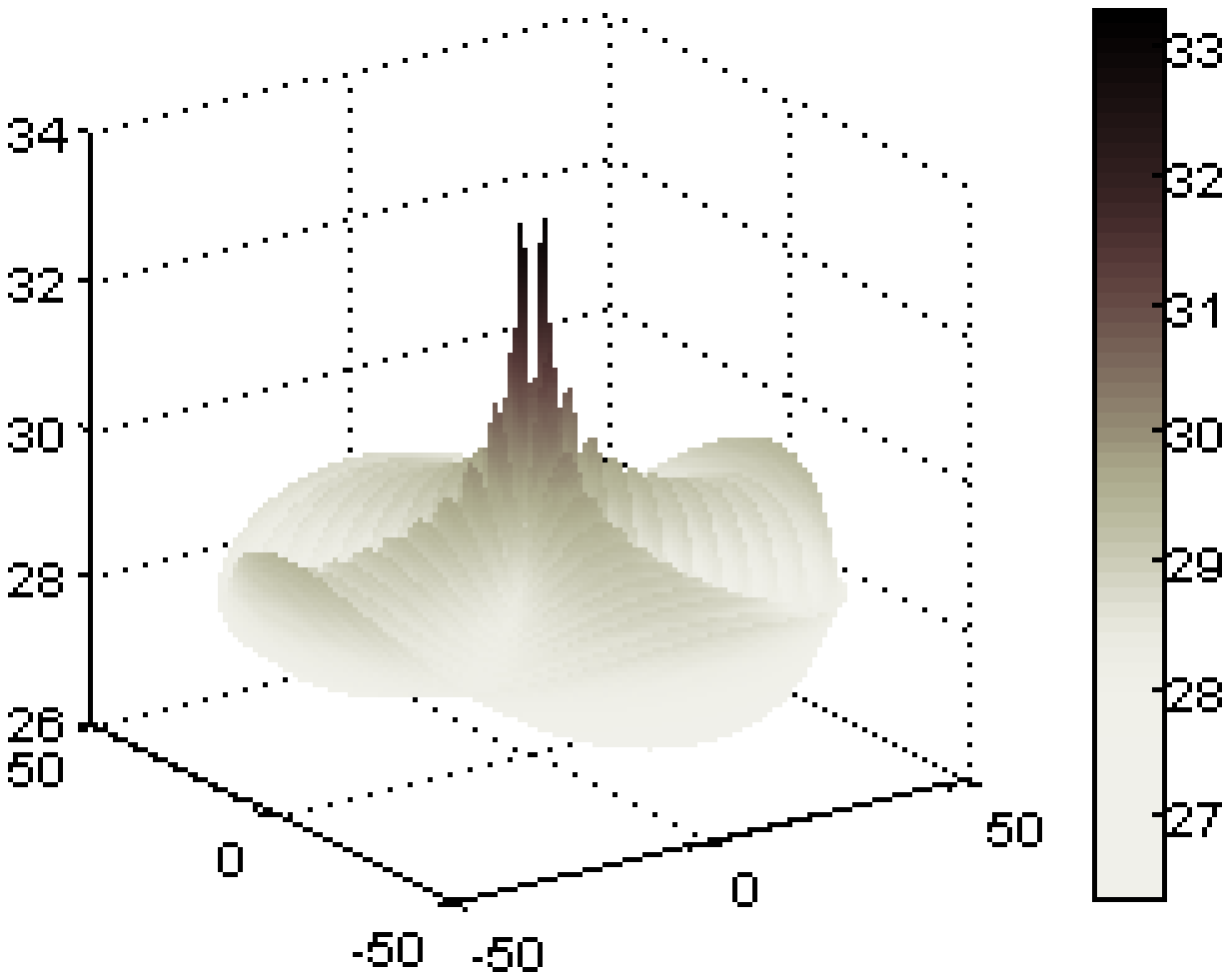} &
\includegraphics[width=0.27\textwidth]{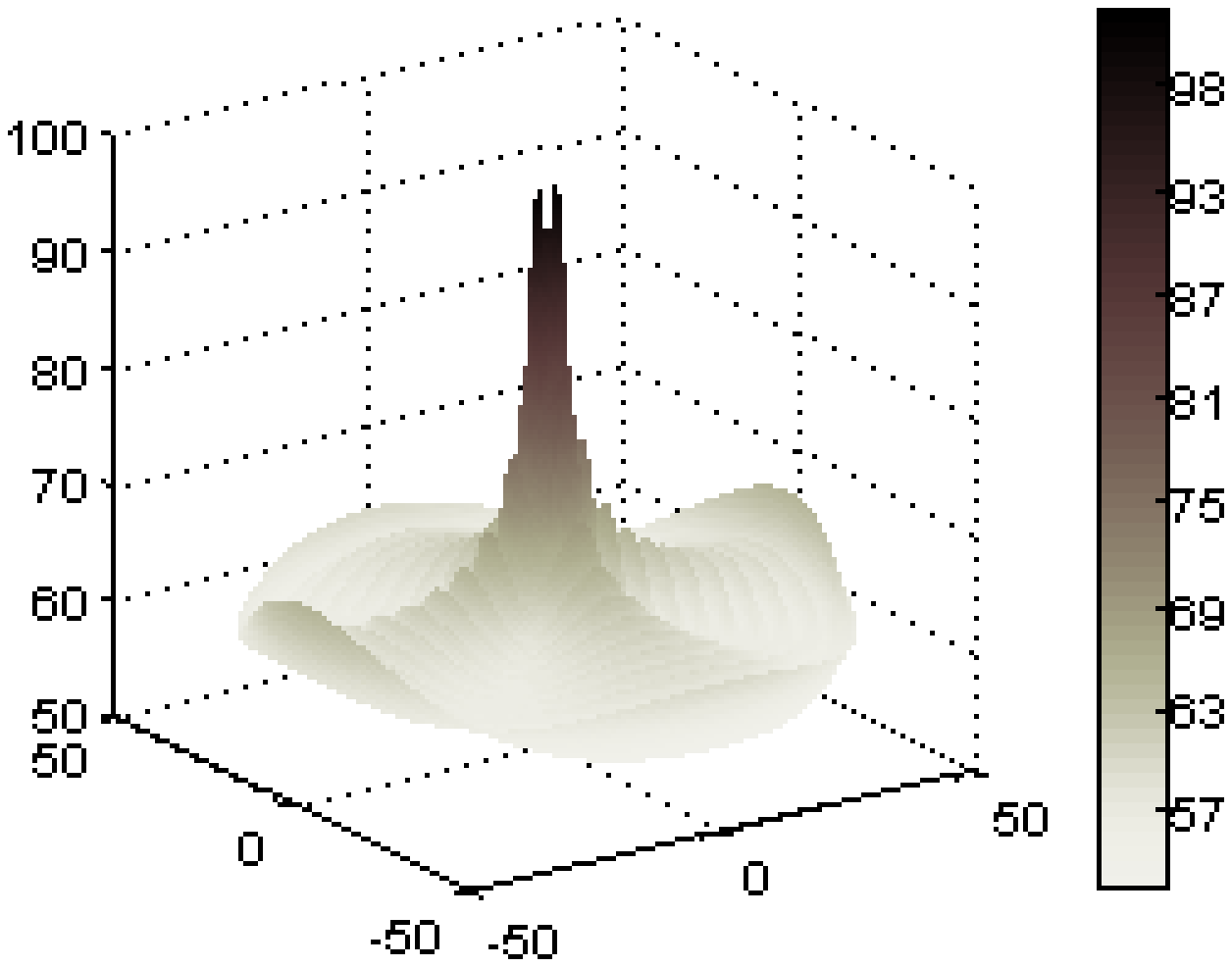} &
\includegraphics[width=0.27\textwidth]{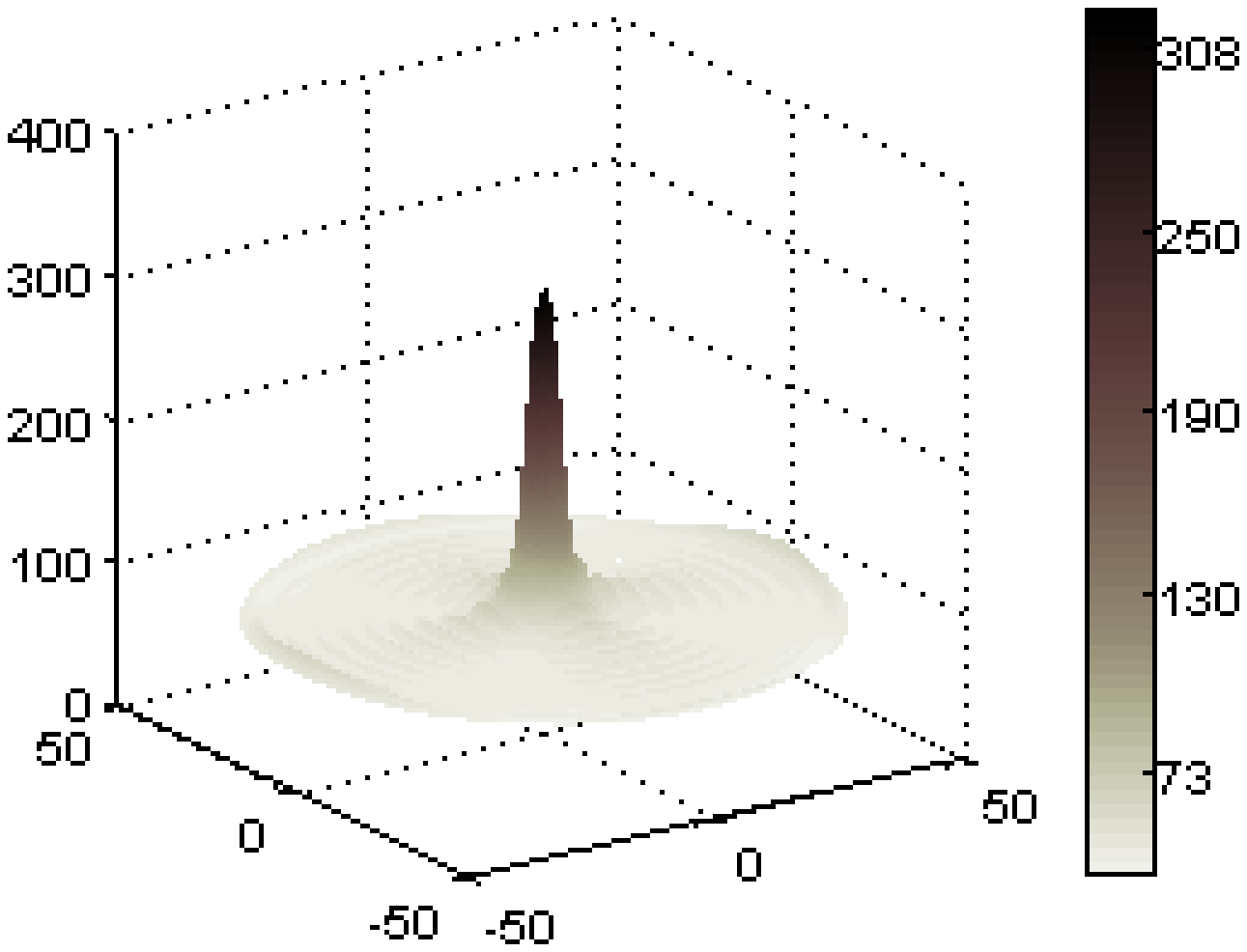} \\
$t=0.499367$ & $t=0.499856$ & $t=0.499914$
\end{tabular}
\end{center}
\caption{\small{Evolution of a complex wave in $S$-regime:
$\sigma=2,\ \beta=3,\ c_0=0,\ k=2,\ T_0= 0.5$.\label{fig14}}}
\end{figure}

Later  the problem of  finding  nonmonotone s.-s.f. in  the $S$-regime,
tending to the nontrivial homogeneous solution $\theta_H$,
 was posed and successfully solved \cite{D26}.
The self-similar equation was solved  with boundary conditions
$$  \ds \frac{\partial\theta_{s,k}}{\partial\xi}=  \ds\frac{1 - \theta_{s,k}}{2\xi} - \gamma
\sqrt{2 \over {\pi \sqrt{\sigma} \xi}}\sin\left(\sqrt{\sigma} \xi - {k\pi\over 2} - {\pi \over 4}\right) \cos(k\phi),
 \ \xi=l\gg 1,  \phi\in\left[0,{2\pi \over k}\right].$$
 The s.-s.f.  for $\beta =\sigma + 1 = 3$ and its evolution in time
are shown in Fig.~\ref{fig14}.

Let us mention, that the existence of continuum of solutions to
the radially symmetric s.-s. problem in
$HS$- and $S$-regimes, which tend to the nontrivial homogenious solution
$\theta_H$ for  $\xi \to \infty$, was mentioned in  \cite{32},
but this result has remained without attention.
As it turned out, it is these solutions that determine the spiral
structures and the complex structures in $HS$- and $S$-regimes.

\section{Open problems}

In spite of the numerous achievements related to the problems considered here,
many interesting questions are still open:
how many complex-symmetry s.-s.f. of the kind $Ei/j$ and $EiMj$,
tending to the trivial constant solution $\theta_s \equiv 0,$ exist;
how can the self-similar problem for the spiral s.-s.f. in the $LS$-regime be closed, and therefore how to
 construct these spiral s.-s.f. numerically;
how wide are the different classes of spiral s.-s.f. for $\beta<\sigma+1$ and the different classes of
complex-symmetry s.-s.f. in $S$- and $HS$-regimes.

Although the method of invariant group analysis
shows that the differential problem admits some kind of self-similar solutions and their numerical realization
is a constructive "proof" of their existence, theoretical proofs are still missing in most of the cases.

The numerical investigation of the accuracy of the approximate solutions on embedded grids shows
optimal-order- and even superconvergence results, but it would be interesting to find theoretical estimates as well.

All these questions pose challenging problems both from theoretical and computational  points of view.

\section*{Acknowledgements}

The  first author is partially  supported by the Sofia University research grant No 181/2012,
the second and the third authors are  partially  supported by the Bulgarian National Science Foundation
under  Grant DDVU02/71.

\end{document}